\pdfoutput=1
\documentclass[runningheads]{llncs}
\usepackage[T1]{fontenc}
\usepackage[english]{babel}
\usepackage{graphicx}
\usepackage{xcolor}
\usepackage{amsmath}
\usepackage{amssymb}
\usepackage{mathabx}
\usepackage{bbm}
\usepackage{etoolbox}
\usepackage{float}
\usepackage{pgfplots}
\usetikzlibrary{graphs,arrows.meta,matrix,fit,patterns,shapes,positioning,fadings}
\pgfplotsset{compat=1.8}
\usepackage{hyperref}
\usepackage[misc,geometry]{ifsym}

\global\newtoggle{QEDHERE}
\AtEndEnvironment{proof}{\iftoggle{QEDHERE}{\global\togglefalse{QEDHERE}}{\phantom{}\qed}}

\newcommand{\midd}[1]{\mathrel{}\middle#1\mathrel{}}

\newcommand{\bbR} {\mathbb{R}}

\newcommand{\bbN} {\mathbb{N}}

\newcommand{\calO} {\mathcal{O}}

\newcommand{\1} {\mathbbm{1}}
\newcommand{\transp} {\scriptscriptstyle \mathsf{T}}

\newcommand{\dx}[1] {\;\mathrm{d}#1}
\newcommand{\mudx}[2] {\;#1\left(\mathrm{d}#2\right)}
\newcommand{\dfx}[2] {\;\mathrm{d}#1\left(#2\right)}

\newcommand{\abs}[1] {\left|#1\right|}
\newcommand{\norm}[1] {\left\lVert #1\right\rVert}
\newcommand{\WD}[2] {\mathrm{WD}\left(#1,#2\right)}

\newcommand{\Prb}[1] {\mathbb{P}\left[#1\right]}

\newcommand{\Ex}[1] {\mathbb{E}\left[#1\right]}
\newcommand{\Exc}[2] {\mathbb{E}_{#1}\left[#2\right]}

\definecolor{tumBlue}{RGB}{0,101,189} % #0065bd
\definecolor{tumDarkBlue}{RGB}{0,82,147} % #005293
\definecolor{tumLightBlue}{RGB}{100,160,200} % #64a0c8
\definecolor{tumLighterBlue}{RGB}{152,198,234} % #98c6ea
\definecolor{tumOrange}{RGB}{227,114,34} % #e37222
\definecolor{tumGreen}{RGB}{162,173,0} % #a2ad00
\definecolor{tumGray}{RGB}{153,153,153} % #999999
\definecolor{tumLight}{RGB}{218,215,203} % #dad7cb

\definecolor{rainbowa}{RGB}{52, 90, 93}
\definecolor{rainbowb}{RGB}{53, 94, 138}
\definecolor{rainbowc}{RGB}{83, 163, 150}
\definecolor{rainbowd}{RGB}{77, 136, 78}
\definecolor{rainbowe}{RGB}{128, 146, 84}
\definecolor{rainbowf}{RGB}{255, 212, 58}
\definecolor{rainbowg}{RGB}{241, 160, 45}
\definecolor{rainbowh}{RGB}{182, 83, 41}
\definecolor{rainbowi}{RGB}{151, 95, 96}
\definecolor{rainbowj}{RGB}{62, 43, 62}

\definecolor{commentpurple}{RGB}{102, 0, 102} % #660066
\colorlet{commentpurplemuted}[RGB]{commentpurple!20!white}

\colorlet{sectionblue}{tumBlue}
\definecolor{linkred}{RGB}{127,0,0} % #7f0000
\definecolor{darklinkred}{RGB}{50,0,0} % #7f0000
\colorlet{headingcolor}{sectionblue}
\colorlet{headingcolormuted}[RGB]{headingcolor!20!white}
\colorlet{linkcolor}{linkred}

\tikzfading[name=fade bottom,
  top color=transparent!0, bottom color=transparent!100, middle color=transparent!50]

\addto\extrasenglish{

}

\begin{document}
\title{Formal Approximations of the Transient Distributions of the M/G/1 Workload Process}
\titlerunning{Transient Distributions of the M/G/1 Workload Process}
\author{Fabian Michel \Letter\inst{1}\orcidID{0009-0005-7768-9111} \and
Markus Siegle\inst{1}\orcidID{0000-0001-7639-2280}}
\authorrunning{F. Michel and M. Siegle}
\institute{Universit\"at der Bundeswehr M\"unchen, Werner-Heisenberg-Weg 39, 85579 Neubiberg, Germany\\
\email{\{fabian.michel,markus.siegle\}@unibw.de}}
\maketitle              % typeset the header of the contribution
\begin{abstract}
  This paper calculates transient distributions of a special class
  of Markov processes with continuous state space and in continuous time,
  up to an explicit error bound. We approximate
  specific queues on $\bbR$ with one-sided Lévy input, such as the M/G/1 workload
  process, with a finite-state Markov chain. The transient distribution
  of the original process is approximated by a distribution with a density
  which is piecewise constant on the state space. Easy-to-calculate
  error bounds for the difference between the approximated and actual
  transient distributions are provided in the Wasserstein distance.
  Our method is fast: to achieve a practically useful error bound,
  it usually requires only a few seconds or at most minutes of computation time.

\end{abstract}

\section{Introduction}

Most of the theory in formal methods for stochastic systems is restricted
to systems where either the state space or the time is discrete. In contrast, we
consider systems where both the state space as well as the time are
continuous. In particular, we would like to calculate transient
distributions of a Markov process in continuous time and with continuous
state space. As exact computations are typically infeasible, we
approximate the transient distributions using a discretization approach,
and we provide formal error bounds for the difference between the actual
and approximated transient distribution. In this paper, we focus on
queues with one-sided Lévy input, as the case of a general Markov process
seems to be very difficult to analyze.

The queues covered by our method include the M/G/1 workload
process. As a motivating example, consider a server setup where the capacity
was chosen such that all jobs can be dealt with reasonably quickly under
the average expected load. Now, assume that a higher than usual job arrival
rate is expected during a short time period (e.g.\ due to new events becoming
available in a booking system, which many users try to book at once). Then,
we can look at the transient workload distribution of the server with
the higher job arrival rate to assess how congested the server will become
in the short heavy load period. Using our error bounds, we can provide
guarantees that the server's workload at a given time will e.g.\ not exceed a given
amount with high probability. In practice, if the
probability of a catastrophic congestion in the heavy load period
is too high, the system administrator could decide to temporarily increase the system's
capacity.

While many results on particular properties of the transient distributions are
available for these types of processes (e.g.\ moments, probability of being idle, see also \autoref{ssec:litrev}),
calculating the transient distributions itself
up to some controllable formal error has received little attention.
Our method can be used to do exactly that and with the transient distribution
available, a variety of questions about the underlying process can be answered.
The error bounds are explicit and easy to calculate, and the computational
cost to meet a predetermined accuracy is reasonably small.
Other approaches, such as the numerical inversion of Laplace transforms
of the transient distributions, usually do not offer error bounds
or only at an unreasonable amount of computational cost. Compared to
the Laplace transform inversion algorithms which are most widely available
in libraries, our method is both faster and more accurate in our numerical
experiments (see \autoref{sec:numex}).

\subsection{Literature review}
\label{ssec:litrev}

Formal error bounds for approximations of Markov processes have been
considered in various settings, but mostly
for models where either space or time are discrete.
An exception is \cite{safetycontspacepurejump}, but it only looks
at pure jump processes.

Next to \cite{safetycontspacepurejump}, one of the most similar papers to the present work is
\cite{precapproxmarkovproc}, which considers a Markov process with general
state space in discrete time. The transition kernel as well as the
initial distribution are assumed to be expressible with a
(Lipschitz-continuous) probability density, the process is approximated
with a finite-state Markov chain and the densities of the
transient distributions are approximated by piecewise constant
densities. This is the same approach that we follow below. However,
we look at continuous-time models and do not assume that the transient
distributions of the original model admit a density. We therefore use a different metric
to measure the error: the Wasserstein distance instead of the $\norm{\cdot}_\infty$-norm
applied to densities as in \cite{precapproxmarkovproc}.
On the other hand, we restrict
ourselves to the state space $\bbR$, unlike \cite{precapproxmarkovproc}.

There are some works, such as \cite{approxmcstochhybrid}, on approximate model checking for stochastic hybrid systems,
which usually have a continuous component in the state space. For more
literature in that direction, also consult the reference lists from \cite{precapproxmarkovproc,approxmcstochhybrid}.
There is also a large body of work on models with discrete state space and continuous time, i.e.\ continuous time
Markov chains. For example, \cite{adaptformalagg,formalbndsstatespaceredmc} provide error bounds for
an approximation of discrete- and continuous-time Markov chains via state
space reduction.

While there seems to be no literature on formal error bounds for approximating
the transient distribution of general Markov processes, some models with
continous time and continuous state space have received considerable attention,
in particular in the analysis of queueing systems.
In \cite{transbehmg1}, some transient characteristics of the M/G/1 workload process
are considered, in particular its moments. \cite{queueslevyfluct} gives a
good overview on the theory behind a more general class of queues with
continuous state space and in continuous time: so-called Lévy-driven queues.
They are defined using a Lévy process (a special type of Markov process
with stationary and independent increments) whose state space is then restricted to
the non-negative reals. This is the setting we focus on, and we will give
more details in \autoref{ssec:levydriven}.

As it turns out, the transient distribution of such queues can often be
characterized by explicit expressions for their
Laplace transforms. Thus,
another approach for calculating transient distributions is the numerical
inversion of these transforms. One should note that the characterizations
are often only given in terms of double or triple transforms, where next to the
Laplace transform of the distribution, additional
transforms in the time variable or in the initial state are considered.
\cite{numtechlevyfluct} follows this approach for the distribution of the running maximum
of a Lévy process, and reports promising results using the inversion
technique from \cite{numtransinvgaussianquad}, albeit without formal error
bounds. \cite{fourierseriesinvtprob,numinvlaplprobdist} propose to use
different inversion techniques and compare the results to estimate the
error, as the computational cost for meeting a pre-defined formal error
bound is often excessive. We will show that a Markov chain approximation
can remedy this problem for transient distributions of Lévy-driven queues.

\subsection{Our contribution}

We present an easy-to-implement method which approximates specific queues with
one-sided Lévy input by a finite-state Markov chain, and which provides explicit
and easy-to-calculate error bounds for the transient distributions in the Wasserstein distance. In particular,
the transient distribution of the original Lévy-driven queue at time points which are multiples of the
discretization parameter $\Delta$ will be approximated by a density which is piecewise
constant on intervals of length $\Delta$. This density is obtained by lifting
the discrete distribution of the approximate model to the original,
continuous state space.
% Using the Wasserstein bound, properties of the
% original process can be proved formally.

\section{Preliminaries}

Consider a general Markov process $X_t$ with a continuous state space
and in continuous time. Assume we want to calculate the transient distribution
at a given point in time. One of the main issues when approximating the
transient distributions is the famous butterfly effect -- small deviations can
result in a completely different future behavior. Therefore, we consider
Lévy-driven queues which offer the advantage that the process behavior
is basically the same everywhere in the state space, making it easier
to control approximation errors.
% One of the main issues when we want to
% approximate the transient distribution and give formal error bounds on
% the approximation is that the (future) behavior of the process can differ arbitrarily strongly even
% for initial states which are very close together. Therefore, when the probability
% mass in our approximation is only slightly shifted, the future
% evolution of the approximated process can be completely off compared to
% the original one. The famous term butterfly effect is often used for
% such a situation. While a generalization of this paper to a wider class
% of processes would be desirable, we will start by considering Lévy-driven
% queues in this work. These offer the advantage that the process behavior
% is basically the same everywhere in the state space, making it easier
% to control approximation errors.

\subsection{Lévy-driven queues}
\label{ssec:levydriven}

We will restrict ourselves to a subclass of Markov processes:
the workload processes of M/G/1 queues, and queues fed
by spectrally negative compound Poisson processes. These belong to
the class of Lévy-driven queues, for which the theory is already well developed. We follow \cite{queueslevyfluct}
to present the most important concepts in this context.

A Lévy process on $\bbR$ is
a Markov process $X_t$ with stationary and independent increments,
and we also require $X_0 = 0$.
Lévy processes can be described by three components:
a deterministic speed $r$ (the process moves with constant speed $r$
upwards or downwards), a Brownian motion
part, and a jump part. We only consider processes without
a Brownian motion part, and which allow only finitely many jumps in a
finite time interval. In addition, we consider so-called spectrally
one-sided Lévy processes, which either only jump upwards or only
jump downwards.
% Our method can be generalized, but we restrict ourselves
% to the above setting for a more concise presentation and leave generalizations
% for future work.

For such Lévy processes (no Brownian motion part, finite jump intensity,
spectrally one-sided), the description is simpler than for general Lévy processes.
The jump rate of $X_t$ into the set $X_t + A$ is defined as $\Pi(A)$,
where $\Pi$ is the so-called Lévy measure on $\bbR \setminus \{0\}$, 
and where $A \subseteq \bbR$ is a Borel set.
The measure $\Pi$ is finite in our case, and either supported on $(0, \infty)$ for upward jumps, or on
$(-\infty,0)$ for downward jumps. The term spectrally positive
process is used for the former and spectrally negative for the latter.

% With Lévy processes defined, we now turn to queues with Lévy input.
Given a spectrally positive or negative Lévy process
with $X_0 = 0$, we define a queue $Q_t$ with net input $X_t$,
started at $Q_0$, by
$Q_t = X_t + \max\left\{Q_0, -\inf_{0\leq s\leq t}X_s\right\}$.
$Q_t$ behaves as a shifted version of $X_t$, except if $Q_t = 0$ and $X_t$
moves down~-- in this case, $Q_t$ stays at $0$.
$Q_t$ is no Lévy process, but it is a Markov process.

$Q_t$ is called the workload of the queue at time $t$, and $X_t$ is called the
net input process, the latter incorporating both the arrivals and the processing of jobs.
% If $X_t$ has a negative drift ($\Ex{X_1} < 0$), then $Q_t$ always
% returns to $0$ and the distribution of $Q_t$ will converge to a stationary
% distribution.
$X_t$ being spectrally positive corresponds to
jobs with varying workloads arriving (according to the measure $\Pi$), and
then being processed at a constant rate $r$ by a server, given by the deterministic
speed of $X_t$. This type of process is also called a compound Poisson process
(see also \cite[page 12, item (2)]{queueslevyfluct}).

We will use the following notation for compound Poisson processes,
both of spectrally positive and negative type: jumps
occur at rate $\lambda = \Pi(\bbR)$, and we let the random variable $\widetilde{B}$
have law $\lambda^{-1} \Pi$. The jump distances are then an iid sequence
with the distribution of the random variable $B = |\widetilde{B}|$.
% Note: $\widetilde{B} \leq 0$ in the spectrally negative case, but the random
% variable $B \geq 0$ just describes the distance that was covered by a jump,
% regardless of the direction of the jump.
In many typical examples, the deterministic speed of the Lévy process
$X_t$ is in the direction opposite to the jump directions.
We thus denote by $r$
the constant speed at which $X_t$ decreases in the spectrally
positive case, while we use $r$ for the speed at which $X_t$
increases in the spectrally negative case.

As mentioned before, the queue arising from a spectrally positive compound
Poisson process can be seen as the workload process of an M/G/1 queue
with server speed $r$, job arrival rate $\lambda$ and job size distribution
given by $B$. The spectrally negative case could for example be used
to model an insurance company which receives premium at rate $r$ per time unit and which has
to pay claims with size distribution $B$, arriving at rate $\lambda$.

Stationary and transient distributions of Lévy-driven queues can be
computed by numerically inverting (single, double, or even triple)
Laplace transforms, as explained in \ref{app:levylaplace}. The discretization approach which we use below has the advantage of
providing better error bounds at a lower computational cost.

\subsection{The Wasserstein distance}

In our approach to approximate the transient distribution of $Q_t$,
we will use the Wasserstein distance to formally bound the distance between the actual transient
distribution and its approximation.
The choice of the Wasserstein distance is deliberate. Other distance
measures such as the total variation distance
often assign the maximal distance to two probability measures which are
orthogonal/singular, which is the case for a Dirac measure and a measure
with a density w.r.t.\ the Lebesgue measure. This would be problematic since we
discretize the original process and approximate its transient distribution
by combinations of uniform distributions over small intervals (see below).
For example, a process started with $Q_0 = x > 0$, i.e.\ a Dirac
measure, would already cause the maximal possible error in the
initial approximation if we used the total variation
distance. Even if the initial distribution is not an issue, jump distributions
with atoms, among others, will not work well in conjunction with such
distances.

The Wasserstein distance is better suited to our approach.
For two probability measures $\mu$ and $\nu$ on $\bbR$ (with the Borel $\sigma$-algebra),
it is defined as:
\begin{align*}
  \WD{\mu}{\nu} = \inf_{\gamma} \int_{\bbR^2} \abs{x - y} \dx{\gamma(x, y)}
  \overset{\textrm{\cite{wassersteindistreals}}}{=} \min_{\gamma} \int_{\bbR^2} \abs{x - y} \dx{\gamma(x, y)}
\end{align*}
where $\gamma$ ranges over all couplings of $\mu$ and $\nu$, i.e.\ we have
$\gamma(A \times \bbR) = \mu(A)$ and $\gamma(\bbR \times A) = \nu(A)$
for measurable $A$ (the marginal distributions of $\gamma$ are $\mu$ and $\nu$, respectively).
The coupling minimizing the above expression describes how to shift
the probability mass of one distribution along the real line in an
optimal way to obtain the second distribution.
We note that, by \cite{wassersteindistreals}, as $\mu$ and $\nu$ are distributions on $\bbR$, this definition is
equivalent to
\begin{align}
  \label{eq:wssrstn_cdf}
  \WD{\mu}{\nu} &= \int_{\bbR} \abs{F_\mu(x) - F_\nu(x)} \dx{x}
\end{align}
where $F_\mu$ and $F_\nu$ are the cumulative distribution functions (CDFs) of
$\mu$ and $\nu$.

% Given the Wasserstein distance between the approximation of
% the transient distribution $\mu$ and the actual transient distribution
% $\nu$, one can give formal bounds for various properties which might
% be of interest. For example, if $f$ is $L$-Lipschitz, then
% we have that $\abs{\Exc{\nu}{f} - \Exc{\mu}{f}} \leq L \cdot \WD{\mu}{\nu}$.
% We also have $\nu\big((a,b)\big) \geq \mu\big((a+\varepsilon, b-\varepsilon)\big) - \varepsilon^{-1} \WD{\mu}{\nu}$
% and $\nu\big((a,b)\big) \leq \mu\big((a-\varepsilon, b+\varepsilon)\big) + \varepsilon^{-1} \WD{\mu}{\nu}$.

\section{Discretization with Formal Error Bounds}
\label{sec:discr}

We will start by approximating the evolution of the workload process $Q_t$ of an
M/G/1 queue with a discrete-time Markov chain. This will allow us to
obtain approximations of the transient distributions of
the process. To simplify notation, we will assume
that the service speed of the M/G/1 queue is fixed at $r = 1$.
As we still allow an arbitrary job arrival rate $\lambda > 0$, this is
no real restriction.

We discretize the model in space
and time, and we truncate the state space to $[0, M]$
with $M > 0$. The precision of the approximation is controlled via the
discretization parameter $\Delta$, and we choose $M$ to be a multiple
of $\Delta$: $M = M_\Delta \cdot \Delta$ with $M_\Delta \in \bbN$.
We approximate $Q_t$ with a discrete-time Markov chain $\widetilde{Q}_k$
on the state space $\{0, 1, \ldots, M_\Delta\}$. The state $\widetilde{Q}_k = 0$
approximates the state $Q_{k\Delta} = 0$ in the original model, while
$\widetilde{Q}_k = i \geq 1$ should hold (approximately) when
$Q_{k\Delta} \in \big((i-1)\Delta, i\Delta\big]$. We discretize
space and time with precisely the same step size $\Delta$ due to the service
speed being $1$. This will be important later.

If $\mu_t$ is the law of $Q_t$, and if $p_k \in \bbR^{M_\Delta + 1}$
is the distribution of $\widetilde{Q}_k$, given by $p_k^{\transp} = p_0^{\transp}P^k$
(with $p_0$ and $P$ still to be defined),
then we approximate $\mu_{k \Delta}$ with
\begin{align}
  \begin{split}
    \widetilde{\mu}_k &:= 
    \Prb{\widetilde{Q}_k = 0} \cdot \delta_0
    + \sum_{i = 1}^{M_\Delta} \Prb{\widetilde{Q}_k = i} \cdot U\big((i-1)\Delta, i\Delta\big) \\
    &= p_k(0) \cdot \delta_0
    + \sum_{i = 1}^{M_\Delta} p_k(i) \cdot U\big((i-1)\Delta, i\Delta\big)
  \end{split}
  \label{eq:approxdist}
\end{align}
where $\delta_0$ is the Dirac measure in $0$ and $U(a,b)$ is a uniform distribution
over the interval $[a, b]$. We later provide a formal bound on
$\WD{\mu_{k \Delta}}{\widetilde{\mu}_k}$, so that $\widetilde{\mu}_k$
(which we can calculate easily) and
this bound can then be used in practice to verify properties of the actual
transient distribution $\mu_{k \Delta}$. Note that $\widetilde{\mu}_k$
is supported on $[0, M]$, while $\mu_{k \Delta}$ is supported on the entire
positive half-line.

The most reasonable choice for $p_0$ is the following:
\begin{align}
  p_0(0) &:= \Prb{Q_0 = 0}, \qquad
  p_0(i) := \Prb{Q_0 \in \big((i-1)\Delta, i\Delta\big]} \;\; (i \geq 1)
  \label{eq:initapprox}
\end{align}
If $M$ is chosen large enough such that the
initial distribution $\mu_0$ of $Q_0$ is supported on $[0, M]$, then this choice
of $p_0$ ensures that $\WD{\mu_0}{\widetilde{\mu}_0} \leq \Delta$: the
probability mass contained in every interval of length $\Delta$ is
correctly represented in $\widetilde{\mu}_0$, and it has to be
shifted by a distance of at most $\Delta$ to obtain $\mu_0$.
Below, we proceed to explain how $P$ should
be calculated, and how to derive error bounds.

\subsection{Transition matrix of discretized M/G/1 model}
\label{ssec:transmat}

We want to choose $P$ such that
\begin{align*}
  P(i,j) \approx \Prb{Q_\Delta \in \big((j-1)\Delta, j\Delta\big] \midd| Q_0 \sim U\big((i-1)\Delta, i\Delta\big)}
\end{align*}
for $i, j \geq 1$ (we have to adapt the expression for the special state $0$).
This ensures that (approximately) the right amount of probability mass is transferred from
the interval $\big((i-1)\Delta, i\Delta\big]$ to the interval $\big((j-1)\Delta, j\Delta\big]$
in the discrete model if the starting distribution is uniform over the
discretization intervals. The distribution of $Q_\Delta$ will
in general not be uniform over these intervals.
This incurs a discretization error at every time step,
as we replace the actual distribution of $Q_\Delta$ with a combination
of uniform distributions in the discretized model.

We can calculate $P(i, j)$ for the M/G/1 queue explicitly, up to a controllable
error. Recall that jobs whose sizes are iid arrive at rate $\lambda$ and are served
at constant speed $r = 1$. We write $F_B$ for the CDF
of $B$, a random variable having the job size distribution.
With probability $e^{-\lambda \Delta}$, no new job arrives within time
$\Delta$ and the probability mass is simply shifted by $\Delta$ downwards
in the state space. With probability $\lambda\Delta e^{-\lambda \Delta}$,
one new job arrives in that same time interval. We will ignore two or more
jobs arriving within the interval $[0, \Delta]$. The reason will become
apparent later~-- basically, it is enough to consider only one job arriving
in order to obtain good error bounds.

\subsubsection{Conditional one-jump CDFs}
Let
\begin{align*}
  F_{\mathrm{oj}}^{(s)}(y) := \Prb{Q_\Delta \leq y \midd| Q_0 = s, 1 \textrm{ job arrival in } [0, \Delta]}
\end{align*}
be the CDF of $Q_\Delta$, conditioned on one jump (oj) in the time
interval $[0, \Delta]$ and started with $Q_0 = s$.
% In contrast to the
% case of no new job arriving, it is a bit more complicated to analyse
% this case.
We further write
\begin{align*}
  F_{\mathrm{oj}}^{[i]}(y) &:= \Prb{Q_\Delta \leq y \midd| Q_0 \sim U\big((i-1)\Delta, i\Delta \big), 1 \textrm{ job arrival in } [0, \Delta]} \quad (i \geq 1) \\
  F_{\mathrm{oj}}^{[0]}(y) &:= \Prb{Q_\Delta \leq y \midd| Q_0 = 0, 1 \textrm{ job arrival in } [0, \Delta]}
  = F_{\mathrm{oj}}^{(0)}(y)
\end{align*}
We will proceed by deriving expressions for $F_{\mathrm{oj}}^{(s)}(y)$
and $F_{\mathrm{oj}}^{[i]}(y)$, which we can then use to calculate $P(i,j)$.
We have to distinguish two cases with respect to $s$.

\subsubsection{Case $s \geq \Delta$}
If $s \geq \Delta$,
then the server will not idle within time $\Delta$ and
\begin{align*}
  F_{\mathrm{oj}}^{(s)}(y)
  &= \Prb{Q_\Delta \leq y \midd| Q_0 = s, 1 \textrm{ job arrival in } [0, \Delta]} \\
  &= \Prb{s + B - \Delta \leq y} = \Prb{B \leq y + \Delta - s} \hspace{2cm} (s \geq 1) \\
  &= F_B(y + \Delta - s)
\end{align*}
This holds because for $Q_\Delta$ to be $\leq y$, we need that the starting
workload $s$ plus the new job size $B$ minus the processed workload 
within the time interval $[0, \Delta]$ (that is, $\Delta$, due to $r = 1$)
is $\leq y$. In consequence,
\begin{align*}
  F_{\mathrm{oj}}^{[i]}(y)
  = \frac{1}{\Delta} \int_{(i-1)\Delta}^{i\Delta} F_{\mathrm{oj}}^{(s)}(y) \dx{s}
  = \frac{1}{\Delta} \int_{y-(i-1)\Delta}^{y-(i-2)\Delta} F_B(s) \dx{s}
  \qquad (i \geq 2)
\end{align*}
Here, we just averaged with respect to the uniform distribution over the
interval $[(i-1)\Delta, i\Delta]$, which is the starting distribution of
$Q_0$ in the definition of $F_{\mathrm{oj}}^{[i]}$.

\subsubsection{Case $s < \Delta$}
For $s < \Delta$, we need to consider that the server might idle some
of the time within the interval $[0, \Delta]$. To simplify calculations,
we will define the idle time as the time spent at $0$ before the new
job arrives (we are still conditioning on one job arrival). It is possible
that $Q_t$ first reaches $0$, then a very small job arrives, and $Q_t$ reaches $0$
again before time $\Delta$. However, it will be easier to consider only
the time spent at $0$ before the arrival as the idle time. In fact, for
the following calculations, we will let $Q_t$ take negative values instead
of being absorbed in $0$, continuing to decrease at constant speed $1$, but \emph{only if} $Q_t$ reaches $0$ \emph{after} the new
job has already arrived. Before the new job arrives, $Q_t$ will be held
at $0$ as before in case the workload $s$ present at time $0$ has already been
processed. $F_{\mathrm{oj}}^{(s)}(y)$ will thus be positive
for $y > -\Delta$, and we still have
\begin{align*}
  F_{\mathrm{oj}}^{(s)}(y) = \Prb{Q_\Delta \leq y \midd| Q_0 = s, 1 \textrm{ job arrival in } [0, \Delta]}
  \qquad (y \geq 0)
\end{align*}
both in the original setting as well as if we let $Q_t$ take negative values
after the job arrival, the equality just doesn't hold for $y < 0$.

The idle time of the server before the new job arrival within $[0, \Delta]$,
if started with workload $s$ at time $0$, is distributed as
$\frac{\Delta - s}{\Delta}U(0,\Delta-s) + \frac{s}{\Delta}\delta_0$:
the time of the new job arrival is distributed uniformly over
$[0, \Delta]$ (when conditioning on one arrival), and thus, with
probability $\frac{s}{\Delta}$, the new job arrives before the old
workload is processed (which would happen at time $s$) and the server does not idle. With
probability $\frac{\Delta - s}{\Delta}$, the job arrives after
$0$ has been reached; then, the idle time is uniformly
distributed between $0$ and $\Delta - s$. The processing time
is distributed as $\Delta$ minus the idle time, i.e.\ its
distribution is $\frac{\Delta - s}{\Delta}U(s,\Delta) + \frac{s}{\Delta}\delta_\Delta$.

We can now write down the equation for $F_{\mathrm{oj}}^{(s)}(y)$:
% considering that the server's processing time is not necessarily $\Delta$ now in the
% case $s < \Delta$:
\begin{align*}
  F_{\mathrm{oj}}^{(s)}(y) &=
  \frac{\Delta - s}{\Delta} \cdot \frac{1}{\Delta - s} \int_s^\Delta F_B(y + t - s) \dx{t}
  + \frac{s}{\Delta} F_B(y + \Delta - s)
  \quad \;\; (y \geq -\Delta) \\
  &= \frac{1}{\Delta} \int_s^\Delta F_B(y + t - s) \dx{t}
  + \frac{s}{\Delta} F_B(y + \Delta - s)
\end{align*}
Note: we averaged over the possible processing times of the server,
and the factor $\frac{1}{\Delta - s}$ in the first line originates
from the density of the distribution $U(s,\Delta)$.

The above expression directly yields $F_{\mathrm{oj}}^{[0]}(y)$:
\begin{align*}
  F_{\mathrm{oj}}^{[0]}(y) = F_{\mathrm{oj}}^{(0)}(y)
  = \frac{1}{\Delta} \int_0^\Delta F_B(y + t) \dx{t}
  = \frac{1}{\Delta} \int_y^{y+\Delta} F_B(s) \dx{s}
\end{align*}
Furthermore, we have
\begin{align*}
  F_{\mathrm{oj}}^{[1]}(y)
  &= \frac{1}{\Delta} \int_0^\Delta F_{\mathrm{oj}}^{(s)}(y) \dx{s} \\
  &= \frac{1}{\Delta^2} \int_0^\Delta \left(\int_s^\Delta F_B(y + t - s) \dx{t}
  + s F_B(y + \Delta - s)\right) \dx{s} \\
  &= \ldots = \frac{2}{\Delta^2} \int_y^{y+\Delta} (y + \Delta - s) F_B(s) \dx{s}
\end{align*}
The final expression can be obtained by exchanging the order of the inner
and the outer integral, as well as by a linear substitution in the integration
variables.

\subsubsection{Calculating $P(i,j)$}
We can use the CDFs from above for a first approximation
\begin{align*}
  \widecheck{P}(i, j)
  := e^{-\lambda\Delta} \bigg(\1_{\{j = i-1 \,\lor\, i=j=0\}}
  + \lambda\Delta \Big(F^{[i]}_{\mathrm{oj}}\big(j\Delta\big) - F^{[i]}_{\mathrm{oj}}\big((j-1)\Delta\big)\Big)\bigg)
\end{align*}
The indicator function corresponds to the case that no jobs arrive
(in which case the probability mass simply shifts one discrete state
to the left), and the second summand to the case with one job arrival~--
more job arrivals are ignored in this approximation. As we ignore more jumps
and as we cut off jumps out of the truncated state space, $\widecheck{P}$
will be a substochastic matrix. We define $P$ by $P = \widecheck{P} + D$
where $D \geq 0$ is a diagonal matrix such that $P$ is stochastic.

We can make the above expression for $\widecheck{P}(i, j)$ more explicit: for $i \geq 2$,
\begin{align*}
  F^{[i]}_{\mathrm{oj}}\big(j\Delta\big) - F^{[i]}_{\mathrm{oj}}\big((j-1)\Delta\big)
  &= \frac{1}{\Delta} \int_{(j-i+1)\Delta}^{(j-i+2)\Delta} F_B(s) \dx{s}
  - \frac{1}{\Delta} \int_{(j-i)\Delta}^{(j-i+1)\Delta} F_B(s) \dx{s} \\
  &= \frac{1}{\Delta} \left(
    \int_{(j-i)\Delta}^{(j-i+1)\Delta} \big(F_B(s+\Delta) - F_B(s)\big) \dx{s}
  \right)
\end{align*}
and hence (equivalent calculations can be done for $i=0,i=1$)
\begin{align*}
  \widecheck{P}(i, j)
  &= e^{-\lambda\Delta} \bigg(\1_{\{j = i-1\}}
  + \lambda \int_{(j-i)\Delta}^{(j-i+1)\Delta} \big(F_B(s+\Delta) - F_B(s)\big) \dx{s}\bigg)
  \quad (i \geq 2) \\
  \widecheck{P}(0, j)
  &= e^{-\lambda\Delta} \bigg(\1_{\{j = 0\}}
  + \lambda \int_{(j-1)\Delta}^{j\Delta} \big(F_B(s+\Delta) - F_B(s)\big) \dx{s}\bigg)\\
  \widecheck{P}(1, j)
  &= e^{-\lambda\Delta} \bigg(\1_{\{j = 0\}}
  + \frac{2\lambda}{\Delta} \int_{(j-1)\Delta}^{j\Delta} (j\Delta - s)\big(F_B(s+\Delta) - F_B(s)\big) \dx{s}\bigg)
\end{align*}
To find $\widecheck{P}$, we thus need to integrate the function
$s \mapsto F_B(s+\Delta) - F_B(s) = \Prb{s < B \leq s + \Delta}$
(for $\widecheck{P}(1, j)$, we actually calculate a convolution
with a piecewise linear triangle function and not just a simple integral).
Depending on the distribution of $B$, we might be able to derive
exact expressions for these integrals, otherwise we use numerical
integration.

% We can also write $\widecheck{P}$ in terms of the density $f_B$ of $B$
% if the distribution of $B$ admits a density, but we will not go into the
% details here.

\subsection{Transition matrix of discretized spectrally negative model}
\label{ssec:transmat_specneg}

Assume now that $Q_t$ is the Lévy-driven queue fed by a spectrally negative
Lévy process $X_t$. $X_t$ is a compound Poisson process with constant
upwards speed $r = 1$ and with downward jumps occurring at rate $\lambda$, the jump
sizes being iid with the distribution of the random variable $B$.

We discretize the state space exactly as in the M/G/1 case, described
at the beginning of \autoref{sec:discr}.
The discretized state $0$ can be dropped in the spectrally negative case,
as $0$ will be left immediately if a jump down to $0$ occurs, due to the
constant positive speed of $1$. However, in some situations, it might make
sense to make the state $0$ absorbing in the spectrally negative case,
corresponding e.g.\ to an insurance company going bankrupt. In such a
case, we would keep the discretized state $0$ (and we would of course also
have to adapt the transition probabilities of the discrete model).

The calculations here are simpler than in the M/G/1 case, and
can be found in \ref{app:transmat_specneg}.
We also end up with a transition matrix $P$
of the discrete model (indexed by indices $1$ through $M_\Delta$,
if we drop state $0$), defined as $P = \widecheck{P} + D$ where $D \geq 0$ is a diagonal
matrix ensuring stochasticity and
\begin{align*}
  \widecheck{P}(i, j) %&= e^{-\lambda \Delta}\left(\1_{\{j = i+1\}} + \lambda\Delta\Big(F^{[i]}_{\mathrm{oj}}\big(j\Delta\big) - F^{[i]}_{\mathrm{oj}}\big((j-1)\Delta\big)\Big)\right)
  %\qquad \quad \;\;\; (i, j \geq 1) \\
  &= \begin{cases}
    \displaystyle e^{-\lambda \Delta}\left(\1_{\{j = i+1\}} + \lambda \int_{(i-j)\Delta}^{(i-j+1)\Delta} \big(F_B(s+\Delta) - F_B(s)\big) \dx{s}\right) & \textrm{ if } j \geq 2 \\
    \displaystyle e^{-\lambda \Delta} \cdot \lambda \left(\Delta - \int_{(i-1)\Delta}^{i\Delta} F_B(s) \dx{s}\right) & \textrm{ if } j = 1
  \end{cases}
\end{align*}

\subsection{Error bounds}
\label{ssec:ebounds}

We now derive an error bound for every step in the discrete model~--
a bound on how much the difference between the actual transient distribution
and the approximated distribution can increase per step in the Wasserstein
distance. Assume that the process starts with initial law $\mu_0$,
i.e.\ $Q_0 \sim \mu_0$. We are given an approximation $\widetilde{\mu}_0$
of $\mu_0$ via the distribution $p_0$ of $\widetilde{Q}_0$ over the
aggregates/intervals as in \eqref{eq:approxdist}.
We do \emph{not} assume that $p_0$ satisfies \eqref{eq:initapprox} because
we want to apply the analysis below to all time steps and not just the initial one.
Instead, we assume that we have a bound $b_0$ on the Wasserstein distance
$\textrm{WD}(\mu_0, \widetilde{\mu}_0)$.

\begin{figure}[htb]
  \begin{center}
    \begin{tikzpicture}[>={Latex[length=1.5mm,width=1.5mm]}]
      \fill[tumOrange!15!white] (-0.9,0.7) rectangle (8.5,-3.5);
      \node[tumOrange] at (3.8,-3.1) {distributions on $[0, \infty)$ respectively $[0,M]$};
      \fill[sectionblue!10!white] (9,0.7) rectangle (11,-3.5);
      \node[sectionblue,align=center] at (10,-2.9) {discrete\\distrib.};
      \node (A) at (0,0) {$Q_0 \sim \mu_0$};
      \node (B) at (0,-2) {$Q_\Delta \sim \mu_\Delta$};
      \node (C) at (4,0) {$Q'_0 \sim \widetilde{\mu}_0$};
      \node (D) at (4,-2) {$Q'_\Delta \sim \textrm{Law}(Q'_\Delta)$};
      \node (E) at (8,0) {$\widetilde{\mu}_0$};
      \node (F) at (8,-2) {$\widetilde{\mu}_1$};
      \node (G) at (10,0) {$\widetilde{Q}_0 \sim p_0$};
      \node (H) at (10,-2) {$\widetilde{Q}_1 \sim p_1$};
      \draw[linkred,<->] (A) -- node[above] {WD $\leq b_0$} (C);
      \draw[tumGreen,<->] (C) -- node[above] {WD $= 0$} (E);
      \draw[linkred,<->] (B) -- node[below] {WD $\leq b_1$} (D);
      \draw[linkred,<->] (D) -- node[below] {WD $\leq b_2$} (F);
      \draw[->,dashed] (G) -- (E);
      \draw[->,dashed] (H) -- (F);
      \node (OD) at (2, -1) {original dynamics};
      \draw[dashed] (0.1, -1) -- (OD);
      \draw[dashed] (3.9, -1) -- (OD);
      \draw[->] (A) -- (B);
      \draw[->] (C) -- (D);
      \draw[->] (G) -- node[left,align=right] {aggregated dynamics\\via $\widetilde{Q}$ and $P$} (H);
    \end{tikzpicture}
    \caption{Bounding the Wasserstein distance}
    \label{fig:wassersteindistevol}
  \end{center}
\end{figure}
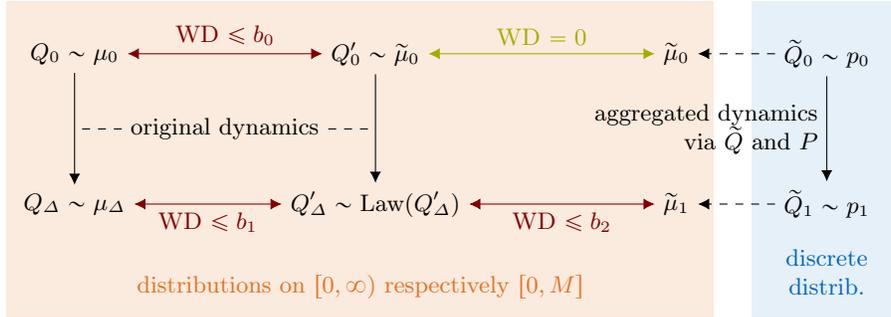

We calculate the distribution
of $\widetilde{Q}_1$ via the matrix $P$, and we want to bound
$\textrm{WD}(\mu_{\Delta}, \widetilde{\mu}_1)$, where $\mu_{\Delta}$ is the
distribution of $Q_{\Delta}$, which we want to approximate with $\widetilde{\mu}_1$,
obtained from the distribution of $\widetilde{Q}_1$. We can apply this bound iteratively to upper bound the Wasserstein
distance $\textrm{WD}(\mu_{k\Delta}, \widetilde{\mu}_k)$ for any $k$ and
therefore give a formal error estimate. We use the strategy depicted in \autoref{fig:wassersteindistevol}:
\begin{itemize}
  \item First, we look at how the error which is already
  present in the initial approximation evolves
  over the time interval $[0, \Delta]$. Consider Markov
  processes $Q$ and $Q'$, started with initial distributions $Q_0 \sim \mu_0$
  and $Q'_0 \sim \widetilde{\mu}_0$, both evolving according to the
  original dynamics of the Lévy-driven queue. Given the
  bound $\textrm{WD}(\mu_0,\widetilde{\mu}_0) \leq b_0$, we will derive
  a bound $b_1$ on $\textrm{WD}(\mu_\Delta,\textrm{Law}(Q'_\Delta))$.
  \item Next, we look at the error caused by approximating
  the dynamics (averaging over the intervals and truncation).
  We will derive a bound $b_2$ on $\textrm{WD}(\textrm{Law}(Q'_\Delta),\widetilde{\mu}_1)$
  where $\widetilde{\mu}_1$ is the distribution as given by \eqref{eq:approxdist} for $k=1$.
  We can calculate the distribution of $\widetilde{Q}_1$ easily
  via the matrix $P$.
  \item By the triangle inequality, we
  can then conclude that $\textrm{WD}(\mu_\Delta, \widetilde{\mu}_1) \leq b_1 + b_2$.
\end{itemize}

\subsubsection{Error caused by initial approximation}
\label{sssec:initerr}

Here, we show that we can choose $b_1 = b_0$ in \autoref{fig:wassersteindistevol}.
That is, if we consider processes started with $Q_0 \sim \mu_0$ and $Q'_0 \sim \widetilde{\mu}_0$,
both evolving according to the same dynamics of the Lévy-driven queue, then
the Wasserstein distance of their transient distributions is bounded by
the initial distance $\textrm{WD}(\mu_0, \widetilde{\mu}_0)$.
We prove this using couplings. We can find a coupling $\gamma$ of $Q_0$ and $Q'_0$
with $\Exc{\gamma}{\abs{Q_0 - Q'_0}} = \textrm{WD}(\mu_0, \widetilde{\mu}_0)$.
We will extend $\gamma$ to a coupling of the two entire processes
(and not just their initial states).

Let $t_1, t_2, \ldots$ be the sequence of times at which the process $Q_t$ (with
$Q_0 \sim \mu_0$) jumps, and let $h_1, h_2, \ldots$ be the corresponding jump
distances. Note: $t_1$ as well as $t_j - t_{j-1}$ for $j \geq 2$ are iid with distribution
$\textrm{Exp}(\lambda)$, independently of the particular value of $Q_0$, and they are in addition independent
of $h_1, h_2, \ldots$. The sequence $h_1, h_2, \ldots$ is itself also an iid
sequence of jump heights with the distribution of $B$.
The jump times and heights of $Q'_t$ (with
$Q'_0 \sim \widetilde{\mu}_0$) follow the same distribution,
for both the M/G/1 and the spectrally negative case.

We can therefore extend the coupling $\gamma$ from the pair $(Q_0, Q'_0)$ to
a coupling $\gamma^\ast$ of the pair
$((Q_t)_{t \geq 0}, (Q'_t)_{t \geq 0})$ by simply letting $Q'_t$ jump with the
same height whenever $Q_t$ jumps. The remaining behavior of both processes is
determined by the constant speed $r = 1$. We look at how
the distance $\abs{Q_t - Q'_t}$ evolves with $t$ under this extended coupling $\gamma^\ast$:
both processes will perform synchronous jumps, and they will both move downwards
with speed $1$ (or upwards in the spectrally negative case), as long as they are not in $0$.
Hence, for the M/G/1 queue, the
distance $\abs{Q_t - Q'_t}$ will stay constant as long as the processes
are either both $> 0$ or both in $0$. When only one process is in $0$,
then the distance $\abs{Q_t - Q'_t}$ decreases with speed $1$, as the other process
will approach $0$ with speed $1$. In the spectrally negative case, the behavior
is similar: $\abs{Q_t - Q'_t}$ will stay constant as long as no jumps occur or if
both processes jump down to a value $> 0$. If one or both processes jump down
to $0$, then $\abs{Q_t - Q'_t}$ goes down abruptly.
Thus, $\abs{Q_t - Q'_t}$ is non-increasing
under the coupling $\gamma^\ast$ in both cases, and in particular $\abs{Q_t - Q'_t} \leq \abs{Q_0 - Q'_0}$.
Therefore (recall that $Q_\Delta \sim \mu_\Delta$),
\begin{align*}
  \textrm{WD}(\mu_\Delta, \textrm{Law}(Q'_\Delta)) &\leq
  \Exc{\gamma^\ast}{\abs{Q_\Delta - Q'_\Delta}}
  \leq \Exc{\gamma^\ast}{\abs{Q_0 - Q'_0}} \\
  &= \Exc{\gamma}{\abs{Q_0 - Q'_0}} = \textrm{WD}(\mu_0, \widetilde{\mu}_0)
\end{align*}
Hence, we can indeed choose $b_1 = b_0$ in \autoref{fig:wassersteindistevol}.
% As a side remark, note that the argumentation above actually extends to all Lévy-driven queues.
% In addition, under the given coupling, we will even have
% $Q_t = Q'_t$ with probability $1$ for some large enough $t$ in the
% setting where the processes always return to $0$, i.e.\ where the net input process
% $X_t$ has a negative drift $\Ex{X_1}$.

\subsubsection{Error caused by aggregated dynamics}

Here, we derive a bound $b_2$ for \autoref{fig:wassersteindistevol}.
That is, we let $Q'_0$ start with distribution $\widetilde{\mu}_0$ as obtained
from a given $p_0$ using \eqref{eq:approxdist} and we then want to bound the
distance between the law of $Q'_\Delta$ (where $Q'_t$ evolves according to
the original process dynamics) and $\widetilde{\mu}_1$ as obtained from $p_1$, where
$p_1^{\transp} = p_0^{\transp} P$ (with $P$ as defined in \autoref{ssec:transmat} or \autoref{ssec:transmat_specneg}).
We thus consider the error caused by approximating the density of $Q'_\Delta$
with a density which is piecewise constant over the aggregation intervals,
and by approximating the transition probabilities between the aggregates
by $P$.

There is no error when $0$ jumps occur in the time interval $[0,\Delta]$,
except in the spectrally negative case for the rightmost discrete state:
if no jump occurs, the probability mass in that state would move out
of the truncated state space. However, we will consider the error caused
by truncation separately below. In principle, if we ignore truncation effects,
the approximation of the density of $Q'_\Delta$ (started with $\widetilde{\mu}_0$ and conditioned on no jump)
is exact.

In contrast, there is an approximation error in the one-jump densities:
the total probability mass in every aggregate is correct,
as we defined $P$ this way, but assuming that it is uniformly distributed
over the intervals is an approximation. The Wasserstein distance of the
piecewise uniform one-jump approximation and the actual distribution of $Q'_\Delta$,
conditioned on one jump, is bounded by $\Delta$, as we only have to redistribute
probability mass within distance $\Delta$ (within one interval) to go
from the approximation to the actual distribution. The probability of
one jump occurring is $\lambda\Delta e^{-\lambda\Delta}$, so the
error per step is at most
\begin{align}
  e_{\textrm{jmpagg}}(\lambda, \Delta)
  = \Delta \cdot \lambda\Delta e^{-\lambda\Delta}
  = \lambda\Delta^2 \cdot \frac{1}{e^{\lambda\Delta}} \leq \lambda\Delta^2
  \label{eq:e_jmpagg}
\end{align}
where the first factor $\Delta$ is the distance by which we have to shift
the probability mass at most to go from one distribution to the other, and where the
second factor $\lambda\Delta e^{-\lambda\Delta}$ corresponds to the amount
of mass we might have to shift. We can further improve the
error bound, as the Wasserstein distance between the piecewise uniform one-jump approximation
and the actual distribution of $Q'_\Delta$, conditioned on one jump, will often be lower than $\Delta$.
See \autoref{sec:numex} for details.

We have a second error source: ignoring
more than one jump per time step of length $\Delta$. The probability
mass moving due to two or more jumps in the original model just stays
where it is in the discretized version. Here, the analysis for the
M/G/1 and the spectrally negative queue differ. We focus on the M/G/1 queue
first. As we allow general
jump height distributions, we might also have to ignore large single
jumps in the M/G/1 case, in particular if arbitrarily large jumps are possible.
The error introduced
by these two types of cut-off can be bounded by:
\begin{gather}
  \begin{gathered}
    \Prb{1 \textrm{ jump in } [0, \Delta]} \cdot \Ex{(\textrm{jump height}) \1_{\{\textrm{jmp.\ hgt.} > M - i\Delta\}} \midd| 1 \textrm{ jump in } [0, \Delta]} \\
    + \sum_{j=2}^\infty \Prb{j \textrm{ jumps in } [0, \Delta]} \cdot \Ex{\textrm{total jump height} \midd| j \textrm{ jumps in } [0, \Delta]}
  \end{gathered}
  \label{eq:e_jmpcut}
\end{gather}
where $i$ is the index of the starting interval in the discrete model.
This follows from the definition of the Wasserstein distance via couplings.
Informally speaking,
we can couple the part of $\textrm{Law}(Q'_\Delta)$ where two or more jumps occurred
in $[0, \Delta]$ or where a single jump led out of the truncated state space with the equal-sized part of $\widetilde{\mu}_1$ resulting
from the amount we added to the diagonal of $P$ to make $\widecheck{P}$ stochastic. The expectation in the
above expression is the integral of the distance of two points w.r.t.\ (a part of) the
coupled measures, as in the definition of the Wasserstein distance. In fact, we could subtract $\Delta$ from the jump distance
within the expectation in most cases due to the constant processing speed $1$.
However, we will not do so as the above expression also gives an upper
bound on the contribution to the Wasserstein distance if the initial distribution
is concentrated on $[0, \Delta]$, where the processing time within time $[0, \Delta]$ is not necessarily $\Delta$.
Rewriting \eqref{eq:e_jmpcut} in terms of $B$, we get (for the M/G/1 case)
\begin{align*}
  &\hphantom{\;=\;} \lambda\Delta e^{-\lambda\Delta} \cdot \Ex{B \1_{\{B > M - i\Delta\}}}
  + \sum_{j=2}^\infty \frac{(\lambda\Delta)^j}{j!}e^{-\lambda\Delta} \cdot j \cdot \Ex{B} \\
  &= \lambda\Delta e^{-\lambda\Delta} \cdot \Ex{B \1_{\{B > M - i\Delta\}}}
  + \lambda\Delta e^{-\lambda\Delta} \cdot \Ex{B} \cdot \sum_{j=2}^\infty \frac{(\lambda\Delta)^{j-1}}{(j-1)!} \\
  &= \underbrace{\lambda\Delta e^{-\lambda\Delta} \cdot \int_{(M - i \Delta,\infty)} x \dfx{F_B}{x}}_{\displaystyle =: e_{\textrm{trunc}}(\lambda, \Delta, i)}
  + \underbrace{\lambda\Delta \left(1 - e^{-\lambda\Delta}\right) \cdot \Ex{B}}_{\displaystyle =: e_{\textrm{jmpcut}}(\lambda, \Delta)} %\\
  % &\leq \lambda\Delta \left(\int_{(M - i \Delta,\infty)} x \dfx{F_B}{x}
  % + \left(1 - e^{-\lambda\Delta}\right) \cdot \Ex{B}\right)
\end{align*}
Note that the Wasserstein error bound only works
if $\Ex{B}$ exists.

For the spectrally negative case, we also ignore two or more jumps per time
interval, but large single jumps are not an issue as jumps cannot go
below $0$. Instead, as previously mentioned, an error occurs when
the probability mass in the topmost discrete space should move out of
the truncated state space due to no jump occurring.
For the error caused by ignoring two or more jumps, we can almost use
the same bound as in the M/G/1 case, but we can take additional advantage
of the fact that jumps are stopped in $0$. As the distribution of $Q'_0$ is supported
on $[0, M]$, no jumps of size larger than $M + \Delta$ can occur
within time $[0, \Delta]$. Therefore, the error caused by ignoring
two or more jumps is bounded by
\begin{align*}
  &\sum_{j=2}^\infty \frac{(\lambda\Delta)^j}{j!}e^{-\lambda\Delta} \cdot \min\{j \cdot \Ex{B}, M + \Delta\} \\
  &\leq \min\left\{\lambda \Delta (1 - e^{-\lambda\Delta}) \Ex{B}, \;\;
  \big(1 - (1 + \lambda\Delta)e^{-\lambda\Delta}\big) (M + \Delta)\right\}
  =: e_{\textrm{jmpcut}}(\lambda, \Delta)
\end{align*}
In fact, we do not need to require that the expectation of $B$ exists in this
case. For the truncation error with respect to the starting interval $i$,
we get
\begin{align*}
  e_{\textrm{trunc}}(\lambda, \Delta, i) = 0 \textrm{~~if } i < M_\Delta,
  \qquad e_{\textrm{trunc}}(\lambda, \Delta, i) = \Delta \cdot e^{-\lambda \Delta} \leq \Delta
  \textrm{~~if }i = M_\Delta
  % \begin{cases}
  %   0 & \textrm{ if } i < M_\Delta \\
  %   \Delta \cdot e^{-\lambda \Delta} \leq \Delta & \textrm{ if } i = M_\Delta
  % \end{cases}
\end{align*}
This is because in the topmost interval (index $M_\Delta$), the mass
which should move upwards by $\Delta$ in case of no jump is $e^{-\lambda\Delta}$
(the probability of no jump).

Putting everything together, we can bound the error per
step in the discrete model by choosing the following $b_2$ in \autoref{fig:wassersteindistevol},
for both the M/G/1 as well as the spectrally negative case (but with
different expressions for $e_{\textrm{jmpcut}}$ and $e_{\textrm{trunc}}$):
\begin{align*}
  b_2 := \sum_{i=0}^{M_\Delta} p_0(i) \cdot \big( e_{\textrm{jmpagg}}(\lambda, \Delta) + e_{\textrm{jmpcut}}(\lambda, \Delta) + e_{\textrm{trunc}}(\lambda, \Delta, i) \big)
\end{align*}
where $p_0$ is the distribution of the discrete model before the current time step.
We want to conclude with an analysis of the behavior of the accumulated error
at time $1$ in the original model (after $\frac{1}{\Delta}$ steps in the discrete model).
For $\Delta \to 0$, the accumulated error should approach $0$ as well,
such that we can actually gain precision by making the aggregation intervals
smaller. If we ignore the truncation part $e_{\textrm{trunc}}(\lambda, \Delta, i)$,
then two remaining parts $e_{\textrm{jmpcut}}(\lambda, \Delta)$
and $e_{\textrm{jmpagg}}(\lambda, \Delta)$ are both
of order $\calO(\Delta^2)$. This is clear for $e_{\textrm{jmpagg}}(\lambda, \Delta)$,
and we have:
\begin{align*}
  e_{\textrm{jmpcut}}(\lambda, \Delta) \leq \lambda\Delta \left(1 - e^{-\lambda\Delta}\right) \cdot \Ex{B}
  = \lambda\Delta \left(\lambda\Delta + \calO(\Delta^2)\right) \cdot \Ex{B}
  = \calO(\Delta^2)
\end{align*}
(For the spectrally negative case, we also have
$e_{\textrm{jmpcut}}(\lambda, \Delta) = \calO(\Delta^2)$ if $\Ex{B}$ does
not exist).
% The error caused by the truncation of single large jumps (M/G/1 case) respectively
% by moving out of the rightmost state (spectrally negative case),
$e_{\textrm{trunc}}(\lambda, \Delta, i)$ is of order
$\calO(\Delta)$ (for fixed $M$). The only requirement for the Wasserstein bound
to be usable in practice is that the error made in the approximation
of the densities of $Q'_t$ per step in the discrete model is
\begin{itemize}
  \item $\calO(\Delta^2)$ for the density approximations conditioned on zero jumps
  (which is true if there is no error in the zero-jump approximation as in our case)
  \item $\calO(\Delta)$ for the density approximations conditioned on one jump
  (which is true if the probability per aggregate is correct in the one-jump approximation as in our case)
\end{itemize}
As a jump only occurs with a probability of $\calO(\Delta)$ within time
$[0, \Delta]$, this implies that the total error per time step
is at most $\calO(\Delta^2)$. This, in turn, implies that the error at
original time $1$ (after $\frac{1}{\Delta}$ steps in the discretized model)
is $\calO(\Delta)$, i.e.\ it does get smaller if we decrease $\Delta$.
This analysis ignores the error due to truncation,
which is a valid approximation in practical settings if the truncation
point is chosen large enough such that only a small part of the probability
mass would have exited the truncated state space within the considered
time horizon. In fact, $e_{\textrm{trunc}}(\lambda, \Delta, i)$ accumulates
to an error of $\calO(1)$ after $\frac{1}{\Delta}$ steps (for $\Delta \to 0$ and $M$ fixed),
but we can make it arbitrarily small by letting $M \to \infty$.

\section{Numerical Example}
\label{sec:numex}

We conclude with a demonstration of the practical applicability
of the presented techniques and error bounds using a numerical example.
% This example also shows that the
% approach using inverse Laplace transforms has problems with discontinuities
% in the probability densities.
% if the most common inversion methods are applied.
% Different distributions for $B$ are considered in the appendix.

The error bounds reported below actually use an improved version of $e_{\textrm{jmpagg}}$
from \eqref{eq:e_jmpagg}: we can calculate the exact CDF of $Q'_\Delta$ in \autoref{fig:wassersteindistevol}
(conditioned on one jump) with the help of the CDFs $F_{\mathrm{oj}}^{[i]}$ obtained in
\autoref{ssec:transmat} (or \ref{app:transmat_specneg} for the spectrally negative case). We can then use
\eqref{eq:wssrstn_cdf} to calculate the Wasserstein distance between
the exact distribution of $Q'_\Delta$ (conditioned on one jump) and the piecewise uniform
approximation, and replace $e_{\textrm{jmpagg}}(\lambda, \Delta)$ by
$\lambda\Delta e^{-\lambda\Delta}$ times the calculated Wasserstein distance.

\begin{figure}[htb]
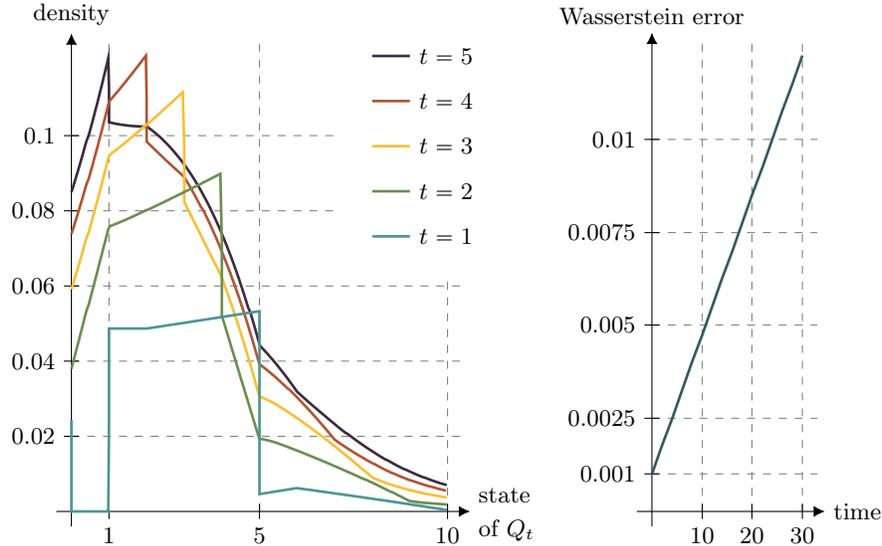

  \begin{center}

    \caption{Transient densities of the M/G/1 workload process started
    with $Q_0 = 1$ at different times $t$. The service speed is
    $1$, the job arrival rate is $\lambda = \frac{1}{4}$, job sizes
    are distributed uniformly over the interval $[1, 5]$.
    Discretization parameter $\Delta = \frac{1}{500}$, truncation
    parameter $M = 50$. On the right: the corresponding Wasserstein
    error bounds.}
    \label{fig:mg1_dirac1_transdens}
  \end{center}
\end{figure}

Consider the M/G/1 queue started at $Q_0 = 1$ with job arrival rate $\lambda = \frac{1}{4}$
and $B$ having a uniform distribution over $[1, 5]$.
This ensures that the process always returns to $0$. \autoref{fig:mg1_dirac1_transdens}
shows how the density of $Q_t$ evolves (the atom at $0$ is not shown).
For example, at time $1$, the density is the sum of the
densities conditioned on a fixed number of jumps, scaled
with the probability of the respective number of jumps ($1 \gg \Delta$, so
our discrete model allows more than one jump up to time $1$). The $1$-jump
part is the uniform distribution over $[1, 5]$ which is very prominent in
for $t = 1$. The $2$-jump part (for $t=1$) has a triangle shape starting at $2$
and going back down to zero at $10$, which is less prominent.
% The $3$-jump part is not visible anymore because
% its amplitude is too small.

The plot on the right of \autoref{fig:mg1_dirac1_transdens} shows how the error bounds
from \autoref{ssec:ebounds} evolve. Here, we used the more
precise version of $e_{\textrm{jmpagg}}$ mentioned above. The initial error $\frac{\Delta}{2} = 0.001$ is
the Wasserstein distance of the Dirac measure at $1$ to a uniform
distribution on the neighboring interval $(1-\Delta, 1]$ of length $\Delta$.
The error increases almost linearly as the truncation error
is comparatively small.

\begin{figure}[htb]
  \begin{center}
    \begin{tikzpicture}[>={Latex[length=1.5mm,width=1.5mm]}]
      \draw[gray,dashed] (0.67, 0) -- (0.67, 3.20);
      \draw[gray,dashed] (3.33, 0) -- (3.33, 3.20);
      \draw[gray,dashed] (6.67, 0) -- (6.67, 3.20);
      \draw[gray,dashed] (10.00, 0) -- (10.00, 3.20);
      \draw[gray,dashed] (0, 1.00) -- (10.20, 1.00);
      \draw[gray,dashed] (0, 2.00) -- (10.20, 2.00);
      \draw[gray,dashed] (0, 3.00) -- (10.20, 3.00);
      \draw[->] (-0.20, 0) -- (10.30, 0) node[right] {state};
      \draw[->] (0, -0.20) -- (0, 3.30) node[above] {density};
      \draw (0.67, 0.1) -- (0.67, -0.1) node[below] {1};
      \draw (3.33, 0.1) -- (3.33, -0.1) node[below] {5};
      \draw (6.67, 0.1) -- (6.67, -0.1) node[below] {10};
      \draw (10.00, 0.1) -- (10.00, -0.1) node[below] {15};
      \draw (0.1, 1.00) -- (-0.1, 1.00) node[left] {0.02};
      \draw (0.1, 2.00) -- (-0.1, 2.00) node[left] {0.04};
      \draw (0.1, 3.00) -- (-0.1, 3.00) node[left] {0.06};
      \draw[rainbowe,line width=1pt] (0.02,0) -- (0.258, -0.0) -- (0.267, 0.0) -- (0.287, -0.0) -- (0.306, -0.0) -- (0.325, -0.0) -- (0.345, -0.0) -- (0.364, -0.0) -- (0.382, -0.0) -- (0.4, -0.0) -- (0.418, 0.0) -- (0.436, 0.001) -- (0.454, -0.001) -- (0.472, -0.002) -- (0.49, 0.006) -- (0.508, -0.008) -- (0.526, 0.006) -- (0.544, -0.003) -- (0.561, -0.001) -- (0.579, 0.004) -- (0.597, -0.002) -- (0.614, -0.011) -- (0.632, 0.051) -- (0.65, -0.165) -- (0.669, 1.494) -- (0.688, 2.642) -- (0.707, 2.343) -- (0.727, 2.46) -- (0.746, 2.449) -- (0.765, 2.4) -- (0.784, 2.462) -- (0.803, 2.425) -- (0.821, 2.423) -- (0.839, 2.447) -- (0.857, 2.433) -- (0.875, 2.426) -- (0.893, 2.438) -- (0.911, 2.437) -- (0.929, 2.43) -- (0.947, 2.432) -- (0.966, 2.436) -- (0.986, 2.434) -- (1.005, 2.432) -- (1.024, 2.434) -- (1.044, 2.435) -- (1.063, 2.434) -- (1.083, 2.433) -- (1.102, 2.434) -- (1.121, 2.434) -- (1.14, 2.434) -- (1.159, 2.434) -- (1.178, 2.434) -- (1.197, 2.434) -- (1.216, 2.434) -- (1.235, 2.434) -- (1.254, 2.434) -- (1.272, 2.434) -- (1.29, 2.434) -- (1.308, 2.434) -- (1.326, 2.434) -- (1.343, 2.435) -- (1.361, 2.437) -- (1.379, 2.439) -- (1.397, 2.441) -- (1.416, 2.443) -- (1.435, 2.445) -- (1.474, 2.45) -- (1.551, 2.459) -- (1.695, 2.475) -- (1.836, 2.491) -- (1.989, 2.509) -- (2.132, 2.525) -- (2.151, 2.527) -- (2.171, 2.529) -- (2.19, 2.531) -- (2.21, 2.533) -- (2.229, 2.537) -- (2.248, 2.54) -- (2.268, 2.542) -- (2.287, 2.542) -- (2.306, 2.543) -- (2.325, 2.545) -- (2.344, 2.549) -- (2.363, 2.554) -- (2.382, 2.557) -- (2.401, 2.559) -- (2.42, 2.557) -- (2.439, 2.556) -- (2.581, 2.569) -- (2.6, 2.567) -- (2.619, 2.57) -- (2.639, 2.581) -- (2.658, 2.595) -- (2.677, 2.605) -- (2.696, 2.606) -- (2.715, 2.598) -- (2.735, 2.584) -- (2.753, 2.574) -- (2.771, 2.573) -- (2.788, 2.584) -- (2.806, 2.604) -- (2.824, 2.626) -- (2.842, 2.641) -- (2.86, 2.644) -- (2.878, 2.631) -- (2.898, 2.606) -- (2.917, 2.58) -- (2.937, 2.568) -- (2.956, 2.577) -- (2.975, 2.608) -- (2.995, 2.651) -- (3.014, 2.688) -- (3.034, 2.704) -- (3.053, 2.689) -- (3.072, 2.647) -- (3.091, 2.591) -- (3.11, 2.543) -- (3.129, 2.525) -- (3.148, 2.551) -- (3.167, 2.622) -- (3.186, 2.72) -- (3.204, 2.809) -- (3.222, 2.866) -- (3.24, 2.859) -- (3.258, 2.765) -- (3.276, 2.575) -- (3.293, 2.293) -- (3.311, 1.939) -- (3.329, 1.542) -- (3.348, 1.108) -- (3.368, 0.714) -- (3.387, 0.397) -- (3.406, 0.179) -- (3.425, 0.063) -- (3.445, 0.037) -- (3.464, 0.078) -- (3.483, 0.155) -- (3.501, 0.235) -- (3.519, 0.303) -- (3.538, 0.347) -- (3.556, 0.363) -- (3.574, 0.353) -- (3.592, 0.324) -- (3.61, 0.287) -- (3.628, 0.25) -- (3.645, 0.223) -- (3.663, 0.209) -- (3.681, 0.21) -- (3.698, 0.223) -- (3.716, 0.245) -- (3.734, 0.271) -- (3.751, 0.295) -- (3.769, 0.315) -- (3.788, 0.327) -- (3.807, 0.329) -- (3.827, 0.322) -- (3.846, 0.31) -- (3.865, 0.297) -- (3.884, 0.284) -- (3.903, 0.276) -- (3.922, 0.274) -- (3.94, 0.278) -- (3.958, 0.285) -- (3.976, 0.295) -- (3.994, 0.306) -- (4.012, 0.315) -- (4.03, 0.322) -- (4.048, 0.325) -- (4.066, 0.325) -- (4.085, 0.32) -- (4.104, 0.313) -- (4.124, 0.304) -- (4.143, 0.294) -- (4.163, 0.286) -- (4.182, 0.28) -- (4.201, 0.276) -- (4.221, 0.275) -- (4.24, 0.276) -- (4.259, 0.278) -- (4.278, 0.281) -- (4.297, 0.283) -- (4.316, 0.284) -- (4.335, 0.284) -- (4.354, 0.283) -- (4.373, 0.28) -- (4.391, 0.276) -- (4.409, 0.272) -- (4.426, 0.267) -- (4.444, 0.263) -- (4.462, 0.259) -- (4.48, 0.255) -- (4.497, 0.253) -- (4.515, 0.251) -- (4.534, 0.25) -- (4.554, 0.249) -- (4.573, 0.249) -- (4.592, 0.249) -- (4.611, 0.248) -- (4.631, 0.248) -- (4.65, 0.246) -- (4.669, 0.245) -- (4.687, 0.242) -- (4.705, 0.24) -- (4.723, 0.237) -- (4.741, 0.234) -- (4.759, 0.231) -- (4.777, 0.228) -- (4.795, 0.226) -- (4.813, 0.223) -- (4.831, 0.221) -- (4.848, 0.22) -- (4.866, 0.218) -- (4.883, 0.216) -- (4.901, 0.215) -- (4.919, 0.214) -- (4.954, 0.211) -- (4.973, 0.21) -- (4.992, 0.208) -- (5.011, 0.206) -- (5.03, 0.204) -- (5.05, 0.202) -- (5.069, 0.199) -- (5.107, 0.195) -- (5.25, 0.178) -- (5.404, 0.162) -- (5.549, 0.146) -- (5.691, 0.13) -- (5.844, 0.113) -- (5.988, 0.097) -- (6.143, 0.079) -- (6.296, 0.062) -- (6.438, 0.045) -- (6.457, 0.043) -- (6.476, 0.04) -- (6.515, 0.035) -- (6.534, 0.033) -- (6.554, 0.031) -- (6.592, 0.027) -- (6.61, 0.025) -- (6.628, 0.023) -- (6.646, 0.022) -- (6.664, 0.021) -- (6.682, 0.02) -- (6.7, 0.019) -- (6.718, 0.018) -- (6.736, 0.018) -- (6.754, 0.017) -- (6.772, 0.017) -- (6.789, 0.017) -- (6.807, 0.017) -- (6.825, 0.017) -- (6.842, 0.017) -- (6.86, 0.017) -- (6.877, 0.017) -- (6.897, 0.017) -- (6.916, 0.017) -- (6.935, 0.017) -- (6.954, 0.017) -- (6.973, 0.017) -- (6.992, 0.017) -- (7.012, 0.017) -- (7.031, 0.017) -- (7.049, 0.016) -- (7.066, 0.016) -- (7.102, 0.016) -- (7.12, 0.016) -- (7.138, 0.015) -- (7.174, 0.015) -- (7.194, 0.015) -- (7.214, 0.014) -- (7.234, 0.014) -- (7.254, 0.014) -- (7.273, 0.014) -- (7.293, 0.014) -- (7.313, 0.014) -- (7.333, 0.014) -- (10, 0);
      \draw[rainbowg,line width=1pt] (0.00, 1.227) -- (0.07, 1.227) -- (0.07, 0.002) -- (0.13, 0.002) -- (0.13, 0) -- (0.60, 0) -- (0.60, 1.176) -- (0.67, 1.176) -- (0.67, 2.437) -- (0.73, 2.437) -- (0.73, 2.441) -- (1.27, 2.441) -- (1.27, 2.442) -- (1.33, 2.442) -- (1.33, 2.448) -- (1.40, 2.448) -- (1.40, 2.454) -- (1.47, 2.454) -- (1.47, 2.461) -- (1.53, 2.461) -- (1.53, 2.468) -- (1.60, 2.468) -- (1.60, 2.475) -- (1.67, 2.475) -- (1.67, 2.482) -- (1.73, 2.482) -- (1.73, 2.489) -- (1.80, 2.489) -- (1.80, 2.496) -- (1.87, 2.496) -- (1.87, 2.502) -- (1.93, 2.502) -- (1.93, 2.509) -- (2.00, 2.509) -- (2.00, 2.516) -- (2.07, 2.516) -- (2.07, 2.523) -- (2.13, 2.523) -- (2.13, 2.530) -- (2.20, 2.530) -- (2.20, 2.537) -- (2.27, 2.537) -- (2.27, 2.544) -- (2.33, 2.544) -- (2.33, 2.551) -- (2.40, 2.551) -- (2.40, 2.558) -- (2.47, 2.558) -- (2.47, 2.565) -- (2.53, 2.565) -- (2.53, 2.571) -- (2.60, 2.571) -- (2.60, 2.578) -- (2.67, 2.578) -- (2.67, 2.585) -- (2.73, 2.585) -- (2.73, 2.592) -- (2.80, 2.592) -- (2.80, 2.599) -- (2.87, 2.599) -- (2.87, 2.606) -- (2.93, 2.606) -- (2.93, 2.613) -- (3.00, 2.613) -- (3.00, 2.621) -- (3.07, 2.621) -- (3.07, 2.628) -- (3.13, 2.628) -- (3.13, 2.635) -- (3.20, 2.635) -- (3.20, 2.642) -- (3.27, 2.642) -- (3.27, 1.472) -- (3.33, 1.472) -- (3.33, 0.219) -- (3.40, 0.219) -- (3.40, 0.222) -- (3.47, 0.222) -- (3.47, 0.229) -- (3.53, 0.229) -- (3.53, 0.237) -- (3.60, 0.237) -- (3.60, 0.244) -- (3.67, 0.244) -- (3.67, 0.251) -- (3.73, 0.251) -- (3.73, 0.258) -- (3.80, 0.258) -- (3.80, 0.265) -- (3.87, 0.265) -- (3.87, 0.272) -- (3.93, 0.272) -- (3.93, 0.276) -- (4.00, 0.276) -- (4.00, 0.273) -- (4.07, 0.273) -- (4.07, 0.267) -- (4.13, 0.267) -- (4.13, 0.260) -- (4.20, 0.260) -- (4.20, 0.254) -- (4.27, 0.254) -- (4.27, 0.247) -- (4.33, 0.247) -- (4.33, 0.241) -- (4.40, 0.241) -- (4.40, 0.234) -- (4.47, 0.234) -- (4.47, 0.228) -- (4.53, 0.228) -- (4.53, 0.221) -- (4.60, 0.221) -- (4.60, 0.215) -- (4.67, 0.215) -- (4.67, 0.209) -- (4.73, 0.209) -- (4.73, 0.202) -- (4.80, 0.202) -- (4.80, 0.196) -- (4.87, 0.196) -- (4.87, 0.189) -- (4.93, 0.189) -- (4.93, 0.183) -- (5.00, 0.183) -- (5.00, 0.176) -- (5.07, 0.176) -- (5.07, 0.170) -- (5.13, 0.170) -- (5.13, 0.163) -- (5.20, 0.163) -- (5.20, 0.157) -- (5.27, 0.157) -- (5.27, 0.150) -- (5.33, 0.150) -- (5.33, 0.143) -- (5.40, 0.143) -- (5.40, 0.137) -- (5.47, 0.137) -- (5.47, 0.130) -- (5.53, 0.130) -- (5.53, 0.123) -- (5.60, 0.123) -- (5.60, 0.117) -- (5.67, 0.117) -- (5.67, 0.110) -- (5.73, 0.110) -- (5.73, 0.103) -- (5.80, 0.103) -- (5.80, 0.096) -- (5.87, 0.096) -- (5.87, 0.089) -- (5.93, 0.089) -- (5.93, 0.083) -- (6.00, 0.083) -- (6.00, 0.076) -- (6.07, 0.076) -- (6.07, 0.069) -- (6.13, 0.069) -- (6.13, 0.062) -- (6.20, 0.062) -- (6.20, 0.055) -- (6.27, 0.055) -- (6.27, 0.048) -- (6.33, 0.048) -- (6.33, 0.041) -- (6.40, 0.041) -- (6.40, 0.034) -- (6.47, 0.034) -- (6.47, 0.027) -- (6.53, 0.027) -- (6.53, 0.020) -- (6.60, 0.020) -- (6.60, 0.015) -- (6.67, 0.015) -- (6.67, 0.013) -- (6.73, 0.013) -- (6.73, 0.012) -- (7.00, 0.012) -- (7.00, 0.011) -- (7.13, 0.011) -- (7.13, 0.010) -- (7.33, 0.010) -- (7.33, 0.009) -- (7.47, 0.009) -- (7.47, 0.008) -- (7.67, 0.008) -- (7.67, 0.007) -- (7.80, 0.007) -- (7.80, 0.006) -- (8.00, 0.006) -- (8.00, 0.005) -- (8.20, 0.005) -- (8.20, 0.004) -- (8.47, 0.004) -- (8.47, 0.003) -- (8.73, 0.003) -- (8.73, 0.002) -- (9.07, 0.002) -- (9.07, 0.001) -- (9.60, 0.001) -- (9.60, 0) -- (10.00, 0);
      \draw[rainbowi,line width=1pt,dashed] (0.00, 1.217) -- (0.02, 0) -- (0.66, 0) -- (0.68, 2.434) -- (1.32, 2.434) -- (1.34, 2.435) -- (1.36, 2.437) -- (1.38, 2.440) -- (1.40, 2.442) -- (1.42, 2.444) -- (1.44, 2.447) -- (1.46, 2.449) -- (1.48, 2.451) -- (1.50, 2.453) -- (1.52, 2.456) -- (1.54, 2.458) -- (1.56, 2.460) -- (1.58, 2.463) -- (1.60, 2.465) -- (1.63, 2.467) -- (1.65, 2.470) -- (1.67, 2.472) -- (1.69, 2.474) -- (1.71, 2.476) -- (1.73, 2.479) -- (1.75, 2.481) -- (1.77, 2.483) -- (1.79, 2.486) -- (1.81, 2.488) -- (1.83, 2.490) -- (1.85, 2.493) -- (1.87, 2.495) -- (1.89, 2.497) -- (1.91, 2.499) -- (1.93, 2.502) -- (1.95, 2.504) -- (1.97, 2.506) -- (1.99, 2.509) -- (2.01, 2.511) -- (2.03, 2.513) -- (2.05, 2.515) -- (2.07, 2.518) -- (2.09, 2.520) -- (2.11, 2.522) -- (2.13, 2.525) -- (2.15, 2.527) -- (2.17, 2.529) -- (2.19, 2.532) -- (2.21, 2.534) -- (2.23, 2.536) -- (2.25, 2.538) -- (2.27, 2.541) -- (2.29, 2.543) -- (2.31, 2.546) -- (2.33, 2.548) -- (2.35, 2.550) -- (2.37, 2.552) -- (2.39, 2.555) -- (2.41, 2.557) -- (2.43, 2.559) -- (2.45, 2.562) -- (2.47, 2.564) -- (2.49, 2.566) -- (2.51, 2.569) -- (2.53, 2.571) -- (2.55, 2.573) -- (2.57, 2.576) -- (2.59, 2.578) -- (2.61, 2.580) -- (2.63, 2.583) -- (2.65, 2.585) -- (2.67, 2.587) -- (2.69, 2.590) -- (2.72, 2.592) -- (2.74, 2.595) -- (2.76, 2.597) -- (2.78, 2.600) -- (2.80, 2.602) -- (2.82, 2.605) -- (2.84, 2.607) -- (2.86, 2.609) -- (2.88, 2.612) -- (2.90, 2.614) -- (2.92, 2.617) -- (2.95, 2.619) -- (2.97, 2.622) -- (2.99, 2.624) -- (3.01, 2.627) -- (3.03, 2.629) -- (3.05, 2.631) -- (3.07, 2.634) -- (3.09, 2.636) -- (3.11, 2.639) -- (3.13, 2.641) -- (3.16, 2.644) -- (3.18, 2.646) -- (3.20, 2.649) -- (3.22, 2.651) -- (3.24, 2.653) -- (3.26, 2.656) -- (3.28, 2.658) -- (3.30, 2.661) -- (3.32, 2.664) -- (3.34, 0.232) -- (3.36, 0.235) -- (3.38, 0.237) -- (3.40, 0.239) -- (3.42, 0.242) -- (3.45, 0.244) -- (3.47, 0.247) -- (3.49, 0.249) -- (3.51, 0.252) -- (3.53, 0.254) -- (3.55, 0.257) -- (3.57, 0.259) -- (3.59, 0.262) -- (3.61, 0.264) -- (3.63, 0.267) -- (3.66, 0.269) -- (3.68, 0.272) -- (3.70, 0.274) -- (3.72, 0.277) -- (3.74, 0.279) -- (3.76, 0.282) -- (3.78, 0.284) -- (3.80, 0.287) -- (3.82, 0.289) -- (3.84, 0.292) -- (3.86, 0.294) -- (3.88, 0.297) -- (3.90, 0.299) -- (3.92, 0.302) -- (3.95, 0.304) -- (3.97, 0.307) -- (3.99, 0.309) -- (4.01, 0.310) -- (4.03, 0.307) -- (4.05, 0.305) -- (4.07, 0.303) -- (4.09, 0.301) -- (4.11, 0.299) -- (4.13, 0.296) -- (4.16, 0.294) -- (4.18, 0.292) -- (4.20, 0.290) -- (4.22, 0.288) -- (4.24, 0.285) -- (4.26, 0.283) -- (4.28, 0.281) -- (4.30, 0.279) -- (4.32, 0.277) -- (4.34, 0.274) -- (4.36, 0.272) -- (4.39, 0.270) -- (4.41, 0.268) -- (4.43, 0.266) -- (4.45, 0.263) -- (4.47, 0.261) -- (4.49, 0.259) -- (4.51, 0.257) -- (4.53, 0.255) -- (4.55, 0.252) -- (4.57, 0.250) -- (4.59, 0.248) -- (4.61, 0.246) -- (4.63, 0.244) -- (4.66, 0.242) -- (4.68, 0.239) -- (4.70, 0.237) -- (4.72, 0.235) -- (4.74, 0.233) -- (4.76, 0.231) -- (4.78, 0.229) -- (4.80, 0.226) -- (4.82, 0.224) -- (4.84, 0.222) -- (4.86, 0.220) -- (4.89, 0.217) -- (4.91, 0.215) -- (4.93, 0.213) -- (4.95, 0.211) -- (4.97, 0.209) -- (4.99, 0.206) -- (5.01, 0.204) -- (5.03, 0.202) -- (5.05, 0.200) -- (5.07, 0.197) -- (5.09, 0.195) -- (5.11, 0.193) -- (5.13, 0.191) -- (5.16, 0.189) -- (5.18, 0.186) -- (5.20, 0.184) -- (5.22, 0.182) -- (5.24, 0.179) -- (5.26, 0.177) -- (5.28, 0.175) -- (5.30, 0.173) -- (5.32, 0.171) -- (5.34, 0.168) -- (5.36, 0.166) -- (5.39, 0.164) -- (5.41, 0.161) -- (5.43, 0.159) -- (5.45, 0.157) -- (5.47, 0.154) -- (5.49, 0.152) -- (5.51, 0.150) -- (5.54, 0.147) -- (5.56, 0.145) -- (5.58, 0.142) -- (5.60, 0.140) -- (5.62, 0.138) -- (5.64, 0.135) -- (5.66, 0.133) -- (5.68, 0.131) -- (5.71, 0.128) -- (5.73, 0.126) -- (5.75, 0.124) -- (5.77, 0.121) -- (5.79, 0.119) -- (5.81, 0.117) -- (5.83, 0.114) -- (5.85, 0.112) -- (5.88, 0.109) -- (5.90, 0.107) -- (5.92, 0.105) -- (5.94, 0.102) -- (5.96, 0.100) -- (5.98, 0.097) -- (6.00, 0.095) -- (6.02, 0.092) -- (6.05, 0.090) -- (6.07, 0.088) -- (6.09, 0.085) -- (6.11, 0.083) -- (6.13, 0.080) -- (6.15, 0.078) -- (6.17, 0.075) -- (6.20, 0.073) -- (6.22, 0.070) -- (6.24, 0.068) -- (6.26, 0.066) -- (6.28, 0.063) -- (6.30, 0.061) -- (6.32, 0.058) -- (6.34, 0.056) -- (6.36, 0.053) -- (6.39, 0.051) -- (6.41, 0.048) -- (6.43, 0.046) -- (6.45, 0.043) -- (6.47, 0.041) -- (6.49, 0.038) -- (6.51, 0.036) -- (6.54, 0.033) -- (6.56, 0.031) -- (6.58, 0.028) -- (6.60, 0.026) -- (6.62, 0.023) -- (6.64, 0.021) -- (6.66, 0.018) -- (6.77, 0.018) -- (6.79, 0.017) -- (6.94, 0.017) -- (6.96, 0.016) -- (7.09, 0.016) -- (7.11, 0.015) -- (7.22, 0.015) -- (7.24, 0.014) -- (7.32, 0.014) -- (7.34, 0.013) -- (7.45, 0.013) -- (7.47, 0.012) -- (7.56, 0.012) -- (7.58, 0.011) -- (7.68, 0.011) -- (7.71, 0.010) -- (7.81, 0.010) -- (7.83, 0.009) -- (7.94, 0.009) -- (7.96, 0.008) -- (8.09, 0.008) -- (8.11, 0.007) -- (8.24, 0.007) -- (8.26, 0.006) -- (8.41, 0.006) -- (8.43, 0.005) -- (8.58, 0.005) -- (8.60, 0.004) -- (8.79, 0.004) -- (8.81, 0.003) -- (9.05, 0.003) -- (9.07, 0.002) -- (9.34, 0.002) -- (9.36, 0.001) -- (9.98, 0.001);
      \fill[white] (7.00, 3.30) rectangle (10.30, 1.10);
      \draw[rainbowe,line width=1pt] (7.50, 2.80) -- (8.00, 2.80) node[right,black] {inverse Laplace};
      \draw[rainbowg,line width=1pt] (7.50, 2.20) -- (8.00, 2.20) node[right,black] {$\Delta = \frac{1}{10}$};
      \draw[rainbowi,line width=1pt,dashed] (7.50, 1.60) -- (8.00, 1.60) node[right,black] {$\Delta = \frac{1}{500}$};
    \end{tikzpicture}
    \caption{Transient densities of the M/G/1 workload process started
    with $Q_0 = 1$ at time $t=1$. The parameters
    are the same as in \autoref{fig:mg1_dirac1_transdens}. Two different
    discretization parameters as well as the inverse Laplace transform approach
    are shown.}
    \label{fig:mg1_dirac1_transdens_comp}
  \end{center}
\end{figure}
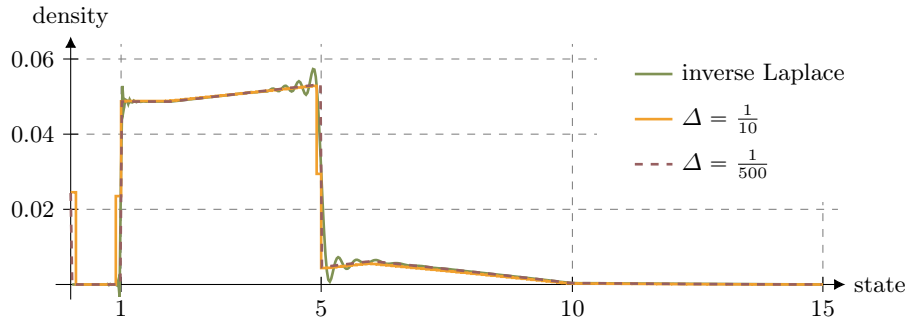

In \autoref{fig:mg1_dirac1_transdens_comp}, we compare setting
$\Delta = \frac{1}{500}$
to $\Delta = \frac{1}{10}$ for $t = 1$.
The density obtained with $\Delta = \frac{1}{10}$ is
is already quite close to the approximation
obtained with $\Delta = \frac{1}{500}$, which shows that even coarse discretizations
can yield good approximations. We also compare with the result obtained
with a double inverse Laplace transform as explained in the appendix. The result obtained by Mathematica \cite{mathematica}
is similar to our results, although without any associated formal error bounds,
and there are oscillatory artifacts near
the discontinuities.

We want to give a short informal account to show that our method
is also attractive with regards to the runtime.
Calculating the transient density approximation with
$\Delta = \frac{1}{10}$ (and the corresponding error bounds) took less than one second on our test machine
(single-threaded, Intel Core i7-1260P CPU at 4.7 GHz),
while Mathematica needs around eight minutes. For $\Delta = \frac{1}{500}$,
the runtime for the discretization approach was around two minutes,
and a common Python library for Laplace transform
inversion, mpmath \cite{mpmath}, did not manage to compute the double inverse at all
in a reasonable amount of time.
Here, a more in-depth comparison, e.g.\ with
the inversion technique from \cite{numtransinvgaussianquad} would be interesting.
% This technique does not offer formal error bounds, though.

The positive density on the interval $(0,0.1]$ for $\Delta = \frac{1}{10}$
in \autoref{fig:mg1_dirac1_transdens_comp}
is a discretization error resulting from ignoring more than one jump per time step.
% when $Q_0 = 1$, then the initial distribution
% is approximated by a uniform distribution on $(0.9,1]$. After time $1$,
% all the mass should have moved to $0$, except if one or more new jobs
% arrived. But as we currently consider job sizes between $1$ and $5$, this
% can only lead to a positive density to the right of $0.9$ (one job of size $1$ arriving
% would mean no change in the workload if one unit of time passes). The
% positive mass in $(0,0.1]$ is a result of ignoring more than one jump per time step of length $\Delta$ and
% just leaving the corresponding probability mass where it was in the approximate
% model.
In \ref{app:numex_heavyloadspecneg}, we also give
an example of an M/G/1 queue under heavy load
and an example of a spectrally negative queue.

\section{Conclusion}

We calculated transient distributions of (a subclass of) queues with
one-sided Lévy input by approximation with a finite Markov chain,
together with explicit error bounds in the Wasserstein distance.
Within a few seconds or minutes of computation time,
the proposed approach can deliver good approximations with error bounds
which are useful in practice. The method is both faster and more accurate
than common inverse Laplace transform approaches, and does not only compute
the transient distribution at a fixed time point as the Laplace approach.
% A more in-depth comparison
% with the Laplace transform approach is still missing.

As a next step, we would like to extend our approach to a wider class
of processes, e.g.\ queues with a Brownian motion part, two-sided
input processes, or an M/G/1 queue with two distinct server
speeds depending on the current load. However, there seem to be fundamental issues when looking
beyond Lévy processes and queues with Lévy input. New methods are required
for these cases.

\begin{credits}
\subsubsection{\discintname}
The authors have no competing interests to declare that are
relevant to the content of this article.
\end{credits}

\appendix
\section*{Appendix}
\setcounter{section}{1}
\renewcommand\thesubsection{Appendix \arabic{subsection}}

\subsection{Lévy-driven queues and Laplace transforms}
\label{app:levylaplace}

Lévy processes
can be characterized by the Laplace transform of the distribution of $X_1$.
For spectrally positive Lévy processes, the Laplace exponent
$\varphi(\alpha) := \ln \Ex{e^{-\alpha X_1}}$ with $\alpha \geq 0$
is well-defined and characterizes the process because
$\Ex{e^{-\alpha X_t}} = e^{t \varphi(\alpha)}$ due to the stationary and
independent increments of a Lévy process. We can write
$\varphi(\alpha) = \alpha r + \int_{(0, \infty)} \left(e^{-\alpha x} - 1 + \alpha x \1_{\{x \in (0,1)\}}\right) \mudx{\Pi}{x}$
where $r$ is the deterministic downwards speed~-- this would be the
server's service speed (rate) in the M/G/1 model. Hence, an explicit expression for $\varphi$
is usually available. Details on the spectrally negative case can be
found in \cite{queueslevyfluct}.

It is possible to derive explicit expressions for the Laplace transform
of the stationary distribution of the queue $Q_t$ in the case $\Ex{X_1} < 0$
where $X_t$ is the net input process (see e.g.\ \cite[Theorem 3.2]{queueslevyfluct}).
% and then to invert the transform numerically to obtain the stationary
% density.
Transient distributions are usually only characterized by double
or triple Laplace transforms, for which the numerical inversion is computationally
much more expensive.
% , in particular if formal error bounds should be met.

\subsubsection{Transient workload of M/G/1 queue}

Consider the following double
transform, where $T_\vartheta$ has an exponential distribution with parameter $\vartheta$
(mean $\vartheta^{-1}$):
\begin{align*}
  \Exc{x}{e^{-\alpha Q_{T_\vartheta}}} = \Ex{e^{-\alpha Q_{T_\vartheta}} \midd| Q_0 = x}
  = \int_0^\infty \vartheta e^{-\vartheta t} \Exc{x}{e^{-\alpha Q_t}} \dx{t}
\end{align*}
By \cite[Theorem 4.1]{queueslevyfluct}, it holds that
\begin{align*}
  \Exc{x}{e^{-\alpha Q_{T_\vartheta}}}
  = \frac{\vartheta}{\vartheta - \varphi(\alpha)}
  \left(e^{-\alpha x} - \frac{\alpha}{\psi(\vartheta)} e^{-\psi(\vartheta) x}\right)
\end{align*}
where $\psi$ is the inverse function of $\alpha \mapsto \varphi(\alpha)$ ($\varphi$ is the Laplace exponent of the net
input process). The existence of $\psi$ is guaranteed in the setting
which we consider here. The density $f_{Q_t}$ of $Q_t$ is thus given by the double inversion
below:
\begin{align*}
  f_{Q_t} = \mathfrak{L}^{-1}\left[\alpha \mapsto \left(
    \mathfrak{L}^{-1}
    \left[ \;\; \vartheta \mapsto \frac{1}{\vartheta - \varphi(\alpha)}
    \left(e^{-\alpha x} - \frac{\alpha}{\psi(\vartheta)} e^{-\psi(\vartheta) x}\right) \;\; \right](t)
  \right)\vphantom{\frac{N^{N^N}}{N^{N^N}}}\right]
\end{align*}
where $\mathfrak{L}^{-1}$ is the inverse Laplace transform.
We first invert $\vartheta \mapsto \frac{1}{\vartheta} \Exc{x}{e^{-\alpha Q_{T_\vartheta}}}$,
the Laplace transform of $t \mapsto \Exc{x}{e^{-\alpha Q_t}}$,
and then invert in $\alpha$ to obtain $f_{Q_t}$.

In the spectrally negative case (see \cite[Section 4.2]{queueslevyfluct}),
the transient distribution can be characterized by a triple transform.
In addition to the Laplace transform of $Q_t$ and the transform
in time, a transform in the initial value is considered.
% Inverting these double or triple transform numerically is possible,
% but computationally expensive, in particular if formal error bounds should be
% met.

\subsection{Transition matrix of discretized spectrally negative model}
\label{app:transmat_specneg}

Here we show how to calculate the matrix $P$ for the discrete
approximation of $Q_t$ in the spectrally negative case.
The queue has a constant upwards
speed $r = 1$, with downward jumps occurring at rate $\lambda$, and with the jump
sizes being iid with the distribution of $B$. Recall
that we measure the positive magnitudes of the jumps, i.e.\ $B \geq 0$,
even though the actual jumps will be downward.
First,
\begin{align*}
  F^{(s)}_{\mathrm{oj}}(y) &:= \Prb{Q_\Delta \leq y \midd| Q_0 = s, 1 \textrm{ jump in }[0, \Delta]} \\
  &= \min\left\{1, \frac{y}{\Delta}\right\} \cdot \Prb{B \geq s + \Delta - y}
\end{align*}
$\Prb{B \geq s + \Delta - y}$ equals $\Prb{s + \Delta - B \leq y}$, i.e.\ the
probability that the starting position plus the deterministic increase $\Delta$
minus the random jump height is $\leq y$. If we want to know whether
$Q_\Delta \leq y$ for some $y \geq \Delta$, then it doesn't matter when the jump
occurs or whether it goes down all the way to $0$, it just needs to be large enough
such that $Q_\Delta$ is at most $y$ in the end, i.e.\ the jump needs to
be at least of size $s + \Delta - y$. If $y < \Delta$, it becomes relevant
when the jump happened: if the jump happens too early in the interval $[0, \Delta]$
(before time $\Delta - y$), then, even if it goes down all the way to $0$,
$Q_\Delta$ will exceed $y$ at the end. Hence, we have to ensure that the
jump happens after time $\Delta - y$ as well as that the jump size is large enough.
As the jump time is distributed uniformly in $[0, \Delta]$ when conditioning
on one jump, the probability of the jump occurring after time
$\Delta - y$ is $\frac{y}{\Delta}$. The required jump size for $Q_\Delta \leq y$
is the same as in the case $y \geq \Delta$, so the total probability
of $Q_\Delta \leq y$ is given by the product of the probability of the jump
occurring late enough and the probability of the jump being far enough.
We then get
\begin{align*}
  F^{[i]}_{\mathrm{oj}}(y) &:= \Prb{Q_\Delta \leq y \midd| Q_0 \sim U\big((i-1)\Delta, i\Delta\big), 1 \textrm{ jump in }[0, \Delta]} \qquad (i \geq 1) \\
  &= \frac{1}{\Delta} \cdot \min\left\{1, \frac{y}{\Delta}\right\} \cdot \int_{(i-1)\Delta}^{i\Delta} \Prb{B \geq s + \Delta - y} \dx{s} \\
  &\overset{\circledast}{=} \min\left\{1, \frac{y}{\Delta}\right\} \cdot \frac{1}{\Delta} \int_{(i-1)\Delta}^{i\Delta} \underbrace{\Prb{B > s + \Delta - y}}_{1 - F_B(s + \Delta - y)} \dx{s} \\
  % &= \min\left\{1, \frac{y}{\Delta}\right\} \cdot \frac{1}{\Delta} \int_{(i-1)\Delta}^{i\Delta} \big(1 - F_B(s + \Delta - y)\big) \dx{s} \\
  &= \min\left\{1, \frac{y}{\Delta}\right\} \left(1 - \frac{1}{\Delta} \int_{i\Delta-y}^{(i+1)\Delta-y} F_B(s) \dx{s} \right)
\end{align*}
where $\circledast$ holds as $\Prb{B \geq s} \neq \Prb{B > s}$ only for
at most countably many $s$, and the set of those $s$ has thus Lebesgue measure $0$.

Finally, similarly to the M/G/1 case, we define the transition matrix $P$
of the discrete model (indexed by indices $1$ through $M_\Delta$, recall
that we drop state $0$) as $P = \widecheck{P} + D$ where $D \geq 0$ is a diagonal
matrix ensuring stochasticity and
\begin{align*}
  \widecheck{P}(i, j) &= e^{-\lambda \Delta}\left(\1_{\{j = i+1\}} + \lambda\Delta\Big(F^{[i]}_{\mathrm{oj}}\big(j\Delta\big) - F^{[i]}_{\mathrm{oj}}\big((j-1)\Delta\big)\Big)\right)
  \qquad \quad \;\;\; (i, j \geq 1) \\
  &= \begin{cases}
    \displaystyle e^{-\lambda \Delta}\left(\1_{\{j = i+1\}} + \lambda \int_{(i-j)\Delta}^{(i-j+1)\Delta} \big(F_B(s+\Delta) - F_B(s)\big) \dx{s}\right) & \textrm{ if } j \geq 2 \\
    \displaystyle e^{-\lambda \Delta} \cdot \lambda \left(\Delta - \int_{(i-1)\Delta}^{i\Delta} F_B(s) \dx{s}\right) & \textrm{ if } j = 1
  \end{cases}
\end{align*}

\subsection{Examples for heavy load and spectrally negative model}
\label{app:numex_heavyloadspecneg}

\begin{figure}[htb]
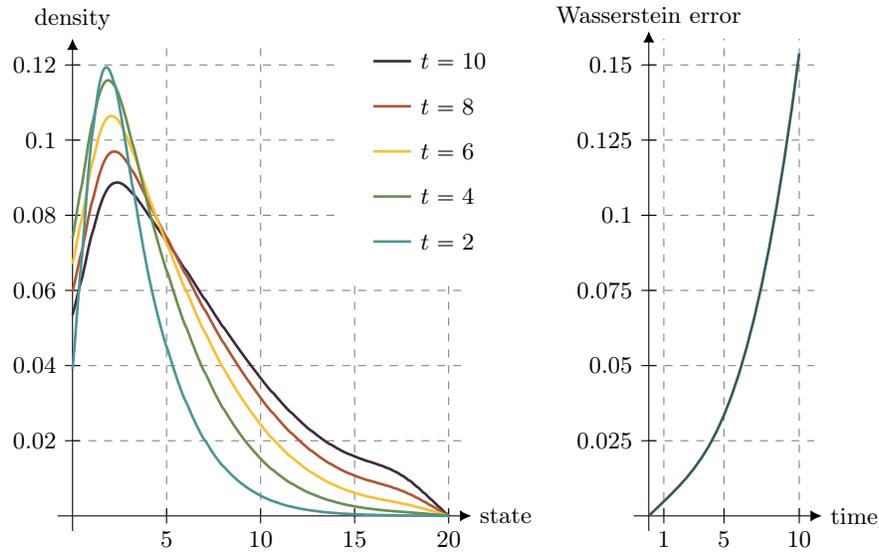

  \begin{center}

    \caption{On the left: transient densities of the M/G/1 workload process 
    at different times $t$, started with $Q_0 = 0$. The constant service speed is
    $1$, the job arrival rate is $\lambda = \frac{2}{5}$, and job sizes
    have an Erlang distribution with scale $6$ and rate $2$ (i.e.\ expectation 3).
    The discretization parameter is $\Delta = \frac{1}{100}$, and the truncation
    parameter is $M = 20$. On the right: the corresponding Wasserstein error
    bounds.}
    \label{fig:mg1_dirac0_erlang_transdens}
  \end{center}
\end{figure}

In \autoref{fig:mg1_dirac0_erlang_transdens}, we look at an example of
the M/G/1 workload process under heavy load. We set
$Q_0 = 0$, we let $B$ have an Erlang distribution with expectation $3$,
and we let $\lambda = \frac{2}{5}$. On average, this results in a workload of
$\frac{6}{5}$ arriving per time unit. As the
server can only process $1$ unit of work per time unit, the workload will
increase to $\infty$ for $t \to \infty$. We see that the formal error
bound in \autoref{fig:mg1_dirac0_erlang_transdens} does increase faster than
linearly. This is due to more and more probability mass accumulating at the
top end of the truncated state space, which causes the truncation error
to increase proportionally. Increasing the truncation point $M$ would result
in an error growth closer to a linear function.

In fact, we can see the effect of the truncation as the densities near
the truncation point $20$ suddenly decrease down to $0$ in
\autoref{fig:mg1_dirac0_erlang_transdens}. This sudden decrease would not
occur in the original model. Still, \autoref{fig:mg1_dirac0_erlang_transdens}
demonstrates that our method can be used to analyze transient workloads
under short periods of heavier-than-usual load. We could easily increase
$M$ (and decrease $\Delta$) to obtain smaller error bounds;
the computation took less than ten seconds here.

\begin{figure}[htb]
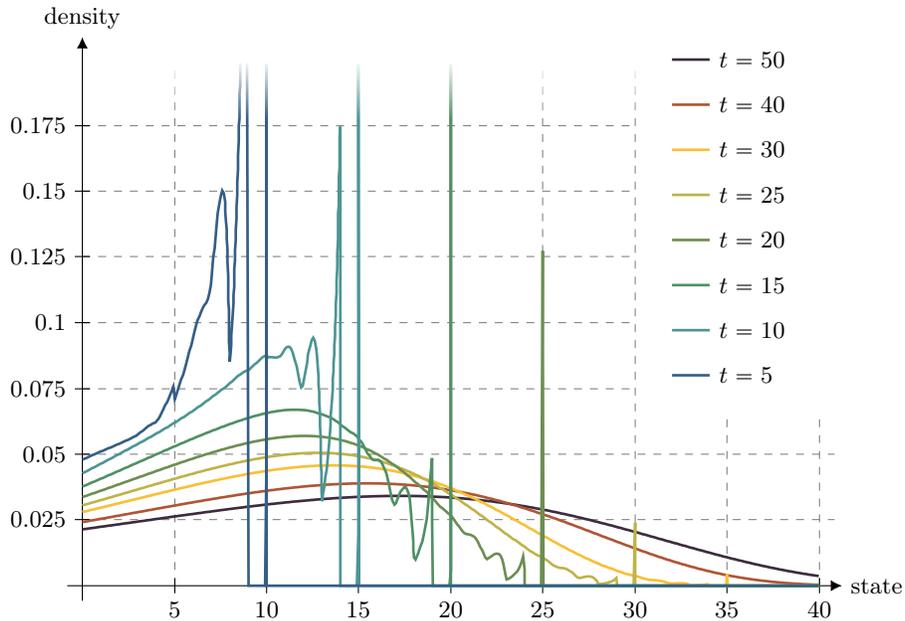

  \begin{center}

    \caption{Transient densities of a Lévy-driven queue started
    with $Q_0 = 5$ at different times $t$. The constant positive speed is
    $1$, the downward jump arrival rate is $\lambda = \frac{1}{3}$, and the jump sizes
    have a Pareto distribution with minimal value $1$ and shape parameter $\alpha = 1.5$.
    The discretization parameter is $\Delta = \frac{1}{100}$, the truncation
    parameter is $M = 55$.}
    \label{fig:insuranceclaims_pareto}
  \end{center}
\end{figure}

In \autoref{fig:insuranceclaims_pareto}, we show an example of a spectrally
negative input process with Pareto-distributed downward jump sizes. The parameters
of the Pareto distribution are chosen such that its expectation is $3$,
and the jump rate is set to $\frac{1}{3}$. As before, the queue moves
with deterministic speed $1$, but now upwards instead of downwards.
The depicted case is thus the critical point
where the average downward jump per time unit is equal to the deterministic
increase per time unit. Here, state $0$ is not absorbing, but one can
see that a significant part of the transient distribution is located
close to $0$. If this process was used to model the capital of an insurance
company, we would conclude that the risk of near-bankruptcy is non-negligible.
The spikes in \autoref{fig:insuranceclaims_pareto} are caused by the
initial Dirac measure moving upwards with speed $1$ (with decreasing
probability mass due to the jumps). In addition, the minimal jump
size of $1$ with the chosen Pareto distribution causes
spikes at integer distances left of the Dirac spike.

\bibliographystyle{splncs04}
\bibliography{references}

\begin{thebibliography}{10}
\providecommand{\url}[1]{\texttt{#1}}
\providecommand{\urlprefix}{URL }
\providecommand{\doi}[1]{https://doi.org/#1}

\bibitem{adaptformalagg}
Abate, A., Andriushchenko, R., \v{C}e\v{s}ka, M., Kwiatkowska, M.: Adaptive formal approximations of {M}arkov chains. Performance Evaluation  \textbf{148}(102207) (2021). \doi{10.1016/j.peva.2021.102207}, \url{https://www.sciencedirect.com/science/article/pii/S0166531621000249}

\bibitem{approxmcstochhybrid}
Abate, A., Katoen, J.P., Lygeros, J., Prandini, M.: Approximate model checking of stochastic hybrid systems. European Journal of Control  \textbf{16}(6),  624--641 (2010). \doi{10.3166/ejc.16.624-641}, \url{https://www.sciencedirect.com/science/article/pii/S0947358010706919}

\bibitem{fourierseriesinvtprob}
Abate, J., Whitt, W.: The {F}ourier-series method for inverting transforms of probability distributions. Queueing Systems  \textbf{10}(1),  5--87 (1992). \doi{10.1007/BF01158520}, \url{https://www.columbia.edu/~ww2040/FourierSeries1992.pdf}

\bibitem{transbehmg1}
Abate, J., Whitt, W.: Transient behavior of the {M/G/1} workload process. Operations Research  \textbf{42}(4),  750--764 (1994). \doi{10.1287/opre.42.4.750}, \url{https://www.columbia.edu/~ww2040/transientworkloadOR94.pdf}

\bibitem{numinvlaplprobdist}
Abate, J., Whitt, W.: Numerical inversion of {L}aplace transforms of probability distributions. ORSA Journal on Computing  \textbf{7}(1),  36--43 (1995). \doi{10.1287/ijoc.7.1.36}, \url{https://www.columbia.edu/~ww2040/LaplaceInversionJoC95.pdf}

\bibitem{numtechlevyfluct}
Asghari, N.M., den Iseger, P., Mandjes, M.: Numerical techniques in {L}\'evy fluctuation theory. Methodology and Computing in Applied Probability  \textbf{16}(1),  31--52 (2014). \doi{10.1007/s11009-012-9296-5}, \url{https://link.springer.com/article/10.1007/s11009-012-9296-5}

\bibitem{queueslevyfluct}
D\k{e}bicki, K., Mandjes, M.: Queues and L\'evy Fluctuation Theory. Springer International Publishing (2015). \doi{10.1007/978-3-319-20693-6_1}, \url{https://link.springer.com/book/10.1007/978-3-319-20693-6}

\bibitem{numtransinvgaussianquad}
den Iseger, P.: Numerical transform inversion using {G}aussian quadrature. Probability in the Engineering and Informational Sciences  \textbf{20}(1),  1--44 (2006). \doi{10.1017/S0269964806060013}, \url{https://papers.ssrn.com/sol3/papers.cfm?abstract_id=1013507}

\bibitem{mpmath}
Johansson, F., et~al.: mpmath: a {P}ython library for arbitrary-precision floating-point arithmetic (version 1.3.0) (March 2023), \url{https://github.com/mpmath/mpmath}

\bibitem{formalbndsstatespaceredmc}
Michel, F., Siegle, M.: Formal error bounds for the state space reduction of {M}arkov chains. Performance Evaluation  \textbf{167},  102464 (2025). \doi{10.1016/j.peva.2024.102464}, \url{https://www.sciencedirect.com/science/article/pii/S0166531624000695}

\bibitem{precapproxmarkovproc}
Soudjani, S.E.Z., Abate, A.: Precise approximations of the probability distribution of a {M}arkov process in time: An application to probabilistic invariance. In: \'Abrah\'am, E., Havelund, K. (eds.) Tools and Algorithms for the Construction and Analysis of Systems, pp. 547--561. Springer (2014). \doi{10.1007/978-3-642-54862-8_45}, \url{https://link.springer.com/chapter/10.1007/978-3-642-54862-8_45}

\bibitem{safetycontspacepurejump}
Soudjani, S.E.Z., Majumdar, R., Abate, A.: Safety verification of continuous-space pure jump {M}arkov processes. In: Chechik, M., Raskin, J.F. (eds.) Tools and Algorithms for the Construction and Analysis of Systems. pp. 147--163. Springer Berlin Heidelberg (2016). \doi{10.1007/978-3-662-49674-9_9}, \url{https://www.cs.ox.ac.uk/people/alessandro.abate/publications/bcSMA16.pdf}

\bibitem{wassersteindistreals}
Vallender, S.S.: Calculation of the {W}asserstein distance between probability distributions on the line. Theory of Probability \& Its Applications  \textbf{18}(4),  784--786 (1974). \doi{10.1137/1118101}

\bibitem{mathematica}
{Wolfram Research, Inc.}: Mathematica, {V}ersion 14.0, \url{https://www.wolfram.com/mathematica}, {C}hampaign, Illinois, 2024

\end{thebibliography}

\end{document}